# Árboles de Forzamiento Semántico Trivalentes para el sistema deductivo paraconsistente P1

## Trivalent Semantic Forcing Trees for Deductive Paraconsistent System P1

**Manuel Sierra-Aristizabal.** msierra@eafit.edu.co

### Resumen |

La semántica de tablas de verdad trivalentes para el sistema de lógica paraconsistente P1, es caracterizada por una herramienta de inferencia visual llamada árboles de forzamiento semántico trivalentes. Dada una fórmula, con esta herramienta se marcan los nodos del árbol asociado a la misma, y se determina si la fórmula es válida o no. En el caso que la fórmula sea inválida, la asignación de valores de verdad que la refuta está determinada por las marcas de las hojas en su árbol de forzamiento.

**Palabras Clave:** árbol de forzamiento, afirmación por defecto, paraconsistencia, semántica, tablas de verdad trivalentes.

### Abstract |

The semantics of trivalent truth tables for the paraconsistent logic system P1, is characterized by a visual inference tool called trivalent semantic forcing trees. Given a formula, with this tool the nodes of the corresponding tree are marked, and it is determined whether the formula is valid or not. In case the formula is invalid, the assignment of truth values that refutes it is determined by the marks of the leaves in its forcing tree.

**Keywords:** forcing tree, truth by default, paraconsistent, semantics, trivalent truth tables.

## 1. INTRODUCCIÓN |

El método de las *tablas semánticas* es presentado en (Beth, 1962), y popularizado por (Smullyan, 1968) con el nombre de *árboles de opciones semánticas*. Con este método se examinan sistemáticamente, todas las opciones que pueden refutar una proposición dada, y se determina si una de estas opciones es lógicamente viable al no generar contradicciones, si esto ocurre se tiene un *contraejemplo* con el cual se refuta la validez de la proposición dada. Si es imposible generar el contraejemplo, es decir, si ninguna de las opciones resulta lógicamente posible, entonces la proposición analizada es válida. Este método ha tenido amplia aceptación, y como hacen (Carnielli, 1987; Barrero y Carnielli, 2005), se ha aplicado a muchos sistemas de lógicas no clásicas.

Los *árboles de forzamiento semántico* presentados por (Sierra, 2001; Sierra, 2006), no buscan construir el contraejemplo explorando todas las opciones posibles, como se hace con las tablas semánticas, sino que, los árboles de forzamiento sólo consideran las opciones que son deductivamente forzadas por las reglas del sistema. En consecuencia, el análisis de validez con los árboles de forzamiento semántico es más sencillo y natural, que el realizado con los árboles de opciones.

En (Sette, 1973; Ciuciura, 2015) el sistema deductivo P1 es presentado, así como su caracterización mediante una semántica trivalente. El sistema P1 es paraconsistente y soporta las contradicciones a nivel de fórmulas atómicas. En la semántica trivalente, el valor de verdad de las fórmulas que no involucran sub-fórmulas atómicas, se determina de la misma forma que se hace con las tablas de verdad del cálculo proposicional clásico, pero los valores de verdad pueden diferir si aparecen fórmulas atómicas, es decir, en la semántica trivalente, el valor de verdad de la negación de una fórmula atómica no necesariamente depende del valor de verdad de la fórmula atómica.

La *semántica de sociedades* biasertivas abiertas es presentada por (Carnielli y Lima-Marques, 1999), a fin de caracterizar el sistema P1. En (Guarín y Montoya, 2003) se extienden los árboles de forzamiento semántico para el cálculo proposicional clásico, con reglas que permiten generar sociedades abiertas, con el fin de refutar fórmulas invalidas del sistema deductivo paraconsistente LPcAt (Lógica Paraconsistente a nivel Atómico) presentado en (Sierra, 2003), y el cual es el resultado de extender el sistema P1 con un operador de incompatibilidad o buen comportamiento de las fórmulas atómicas. En (Sierra, 2018) se presentan los *Árboles de forzamiento semántico para sociedades abiertas*, con las cuales se caracterizan los sistemas P1 y LPcAt.

En este trabajo se presentan los *Árboles de forzamiento semántico trivalentes* para el sistema deductivo P1. Se prueba con todo detalle la equivalencia entre la presentación con árboles de forzamiento trivalentes y la presentación semántica con las tablas de verdad del sistema P1. Finalmente, se muestra que, si una fórmula es inválida, lo cual implica que el árbol de forzamiento de la fórmula está bien marcado, entonces la lectura de las marcas de las fórmulas atómicas proporciona una asignación de valores de verdad, con la cual se refuta la validez de la fórmula.

## 2. LENGUAJE DE LPcAt y P1 |

El lenguaje de los sistemas LPcAt y P1 (excluyendo I), consta de los conectivos binarios: →, ∧ y ∨ (condicional, conjunción y disyunción), de los conectivos monádicos: I, ~ y ¬ (incompatibilidad, negación fuerte y cuestionamiento o negación débil), además del paréntesis izquierdo y el paréntesis derecho. También se tiene una cantidad enumerable de *fórmulas atómicas*. El conjunto de *fórmulas* de LPcAt y P1 (excluyendo R2) es generado por las siguientes reglas y sólo por ellas:

Toda *fórmula atómica* es una *fórmula*.

R1. Si X es una *fórmula* entonces ~(X) y ¬(X) son *fórmulas*.

R2. Si p es una fórmula atómica entonces Ip es una fórmula.

R3. Si X y Y son *fórmulas* entonces (X)∧(Y), (X)∨(Y) y (X)→(Y) son *fórmulas*.

## 3. ÁRBOL DE UNA FÓRMULA |

Sea X una fórmula, el *árbol de* X se representa por *Ar[X]* y se construye utilizando las reglas presentadas en la figura 1, donde X y Y son fórmulas arbitrarias, p es una fórmula atómica.

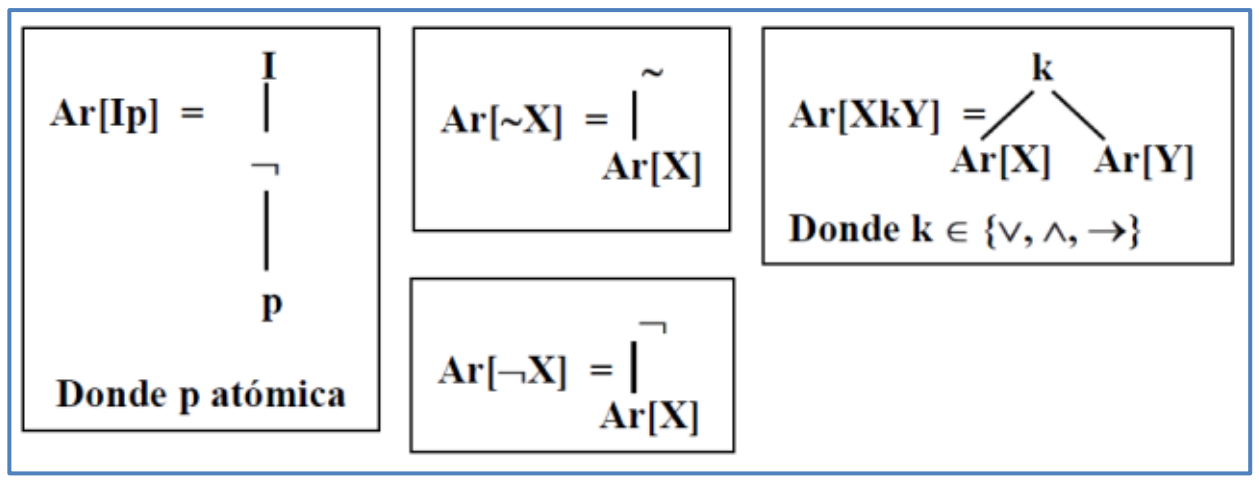

Figura 1

Se define el *árbol del argumento* 'de $X_1$, ... , $X_n$ se infiere Y' como: Ar[($X_1$∧ ... ∧$X_n$)→Y]. El nodo superior del árbol de la fórmula X, es llamado la *raíz* del árbol, se denota R[X] y corresponde al operador principal de la fórmula X. Los nodos inferiores, es decir aquellos de los cuales no salen ramas, son llamados *hojas* y corresponden a las fórmulas atómicas.

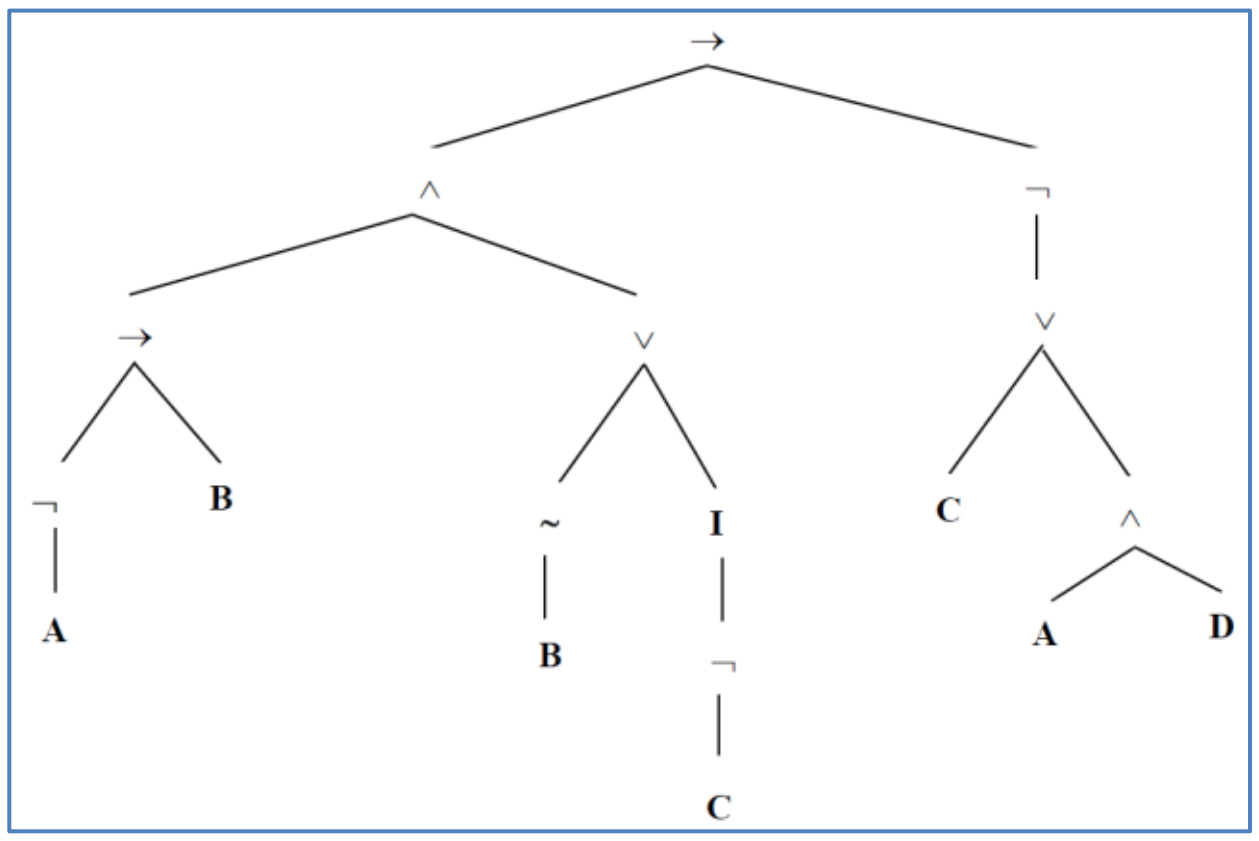

Figura 2

Por ejemplo, para el argumento: 'de ¬A→B y ~B∨IC se infiere ¬[C∨(A∧D)]', donde A, B, C y D son fórmulas atómicas, el condicional asociado es: [(¬A→B) ∧ (~B∨IC)] → ¬[C∨(A∧D)] y su árbol es presentado en la figura 2.

## 4. MARCANDO LOS NODOS DE UN ÁRBOL |

Si un nodo *C* es uno de los conectivos monádicos, I, ~ o ¬, entonces su único hijo se llama el *alcance del operador* y para hacer referencia a él se utiliza la notación *aC*.

Si un nodo *K* es uno de los conectivos binarios ∧, ∨ o →, entonces para sus *hijos izquierdo y derecho* se utiliza la notación *iK* y *dK* respectivamente.

Para toda fórmula Y, el *nodo asociado* a Y es la raíz de Y, *R[Y]*, la cual a su vez es el operador principal de Y en el caso que, Y no sea atómica, o es la misma Y en el caso que Y sea atómica.

Para una fórmula X, H(X) el conjunto de hojas del Ar[X], y N(X) el conjunto de nodos de Ar[X].

Para cada fórmula X, una *función de marca de hojas m* (o simplemente función de marca), es una función de H(X) en {0, *, 1}.

Si *m*(p) = 1 entonces se dice que la hoja p está marcada con 1, o que *es aceptada*.

Si *m*(p) = 0 entonces se dice que la hoja p está marcada con 0, o que *es rechazada*.

Si *m*(p) = * entonces se dice que la hoja p está marcada con *, o que *es aceptada por defecto o por falta de evidencia de lo contrario*.

Cada función de marca de hojas *m*, puede ser extendida de manera única (la prueba se presenta en la proposición 5.13), a una *función de marca de nodos, M*, de N(X) en {0, 1, *}, haciendo *M*(h) = *m*(h) si h es una hoja, y aplicando las reglas primitivas y derivadas para el forzamiento de marca, las cuales son presentadas a continuación, en las secciones 4.1 a 4.7 y 5.1 a 5.12:

*Reglas para marcar los nodos*

$$M(p) = m(p) \text{ si p es una hoja.}$$

Si el nodo *n* no es una hoja entonces M(*n*) solo toma valores en {0, 1}, nunca toma el valor *.

### 4.1 Reglas primitivas para la negación |

A¬. *Aceptación del Cuestionamiento*: Una fórmula no es aceptada cuando es cuestionada.

$$M(\neg) = 1 \Rightarrow M(a\neg) \neq 1$$

Esta regla se ilustra en el paso 8 de la fig 21.

Ra¬. *Rechazo del Alcance del Cuestionamiento*: Las fórmulas rechazadas son cuestionadas.

$$M(a\neg) = 0 \Rightarrow M(\neg) = 1$$

Esta regla se ilustra en el paso 3 de la fig 26.

A*a¬. *Aceptación por defecto del Alcance del Cuestionamiento*: Las fórmulas que son aceptadas por defecto son cuestionadas.

$$M(a\neg) = * \Rightarrow M(\neg) = 1$$

Esta regla se ilustra en el paso 5 de la fig 25.

### 4.2 Reglas primitivas para la incompatibilidad|

RI. *Rechazo de la Incompatibilidad*: Los enunciados atómicos compatibles con su cuestionamiento se aceptan por defecto.

$$M(I) = 0 \Rightarrow M(a\neg) = *$$

Esta regla se ilustra en el paso 4 de la fig 25.

A*a¬I. *Aceptación por defecto del alcance del Cuestionamiento en la Incompatibilidad*: Los enunciados atómicos que son aceptados por defecto son compatibles con su cuestionamiento.

$$M(a\neg) = * \Rightarrow M(I) = 0$$

Esta regla se ilustra en el paso 8 de la fig 27.

### 4.3 Reglas primitivas para el condicional |

R→. *Rechazo del condicional*: Si un condicional es rechazado entonces el consecuente es rechazado y el antecedente no es rechazado.

$$M(\rightarrow) = 0 \Rightarrow M(i\rightarrow) \neq 0 \text{ y } M(d\rightarrow) = 0$$

Esta regla se ilustra en los pasos 2 y 3 de la fig 21 y también en los pasos 6 y 7.

nRiRd→. No *Rechazo a la izquierda y Rechazo a la derecha en el condicional*: Si el antecedente de un condicional no es rechazado y el consecuente es rechazado entonces el condicional es rechazado.

$$M(i\rightarrow) \neq 0 \text{ y } M(d\rightarrow) = 0 \Rightarrow M(\rightarrow) = 0$$

## 4.4 Reglas primitivas para la conjunción |

nRinRd∧. No *Rechazo a la izquierda y no Rechazo a la derecha en la conjunción*: Si ambos coyuntos no son rechazados entonces la conjunción es aceptada.

$$M(i\wedge) \neq 0 \text{ y } M(d\wedge) \neq 0 \Rightarrow M(\wedge) = 1$$

Esta regla se ilustra en el paso 8 de la fig 23.

A∧. Aceptación de *la conjunción*: Si una conjunción es aceptada entonces ambos coyuntos no son rechazados.

$$M(\wedge) = 1 \Rightarrow M(i\wedge) \neq 0 \text{ y } M(d\wedge) \neq 0$$

Esta regla se ilustra en los pasos 7 y 8 de la fig 24.

## 4.5 Reglas primitivas para la disyunción |

RiRd∨. *Rechazo a la izquierda y Rechazo a la derecha en la disyunción*: Si ambos disyuntos son rechazados entonces la disyunción es rechazada.

$$M(i\vee) = 0 \text{ y } M(d\vee) = 0 \Rightarrow M(\vee) = 0$$

R∨. Rechazo de *la disyunción*: Si una disyunción es rechazada entonces ambos disyuntos son rechazados.

$$M(\vee) = 0 \Rightarrow M(i\vee) = 0 \text{ y } M(d\vee) = 0$$

Esta regla se ilustra en los pasos 2 y 3 de la figura 23.

## 4.6 Reglas primitivas para la negación fuerte |

A~. *Aceptación de la Negación fuerte*: Una fórmula es rechazada cuando se acepta su negación.

$$M(\sim) = 1 \Rightarrow M(a\sim) = 0$$

Esta regla se ilustra en el paso 6 de la fig 28.

Ra~. *Rechazo del alcance de la Negación fuerte*: Una negación es aceptada cuando se rechaza su alcance.

$$M(a\sim) = 0 \Rightarrow M(\sim) = 1$$

## 4.7 Validez |

Se dice que una fórmula X es *A-válida* (válida desde el punto de vista de los árboles) si y solamente si para toda función de marca *m*, se tiene que *M*(R[X]) ≠ 0.

Se dice que una fórmula X es *A-inválida* si no es A-válida, es decir si existe una función de marca *m*, tal que *M*(R[X]) = 0. En este caso se dice que la *función de marca refuta* la fórmula X. También se dice que el *árbol de X está bien marcado* (*ABM*, todos sus nodos están marcados de acuerdo a las reglas sin generar contradicciones).

Para el análisis de validez se tienen 2 métodos:

1. **Forzamiento directo**. El objetivo es llegar, aplicando las reglas, a que la raíz está marcada con 1. En este caso se concluye que la fórmula es A-válida, y se justifica como RM1.

2. **Forzamiento indirecto**. En este caso se parte del supuesto que la fórmula es A-inválida, para ello se marca la raíz con 0 y luego se aplican las reglas.

Si en algún momento resulta una contradicción (doble marca) entonces se concluye que la fórmula es A-válida, y se justifica como árbol mal marcado (AMM).

Si se marcan todos los nodos sin contradecir las reglas y no hay doble marca entonces se concluye que la fórmula es A-inválida, y se justifica como árbol bien marcado (ABM).

# 5. REGLAS DERIVADAS PARA EL FORZAMIENTO DE MARCAS |

Las reglas primitivas para el forzamiento de marcas, son suficientes para estudiar las propiedades de los árboles de forzamiento, pero en la práctica, cuando se trata de marcar todos los nodos de un árbol, es importante tener reglas que cubran todas las posibilidades. A continuación, se presenta un juego completo de *reglas derivadas*.

## 5.1. Proposición. *Reglas derivadas para el cuestionamiento* |

Aa¬. *Aceptación del Alcance del Cuestionamiento*: Las fórmulas que son aceptadas no son cuestionadas.

$$M(a\neg) = 1 \Rightarrow M(\neg) = 0$$

Esta regla se ilustra en el paso 13 de la fig 24 y también en el paso 13 de la fig 25.

R¬. *Rechazo del Cuestionamiento*: Si una fórmula no es cuestionada, entonces es aceptada.

$$M(\neg) = 0 \Rightarrow M(a\neg) = 1$$

Esta regla se ilustra en el paso 8 de la fig 22.

Prueba de Aa¬: Sea $M(a\neg) = 1$. Supóngase que $M(\neg) = 1$, entonces por A¬ se infiere $M(a\neg) \neq 1$, lo cual no es el caso. Por lo tanto, forzosamente $M(\neg) = 0$.

Prueba de R¬: Sea $M(\neg) = 0$. Supóngase que $M(a\neg) \neq 1$, entonces por Ra¬ y A*a¬ se infiere $M(\neg) = 1$, lo cual no es el caso. Por lo tanto, forzosamente $M(a\neg) = 1$.

## 5.2. Proposición. *Reglas derivadas para la incompatibilidad |*

A¬AI. *Aceptación del Cuestionamiento y Aceptación de la Incompatibilidad*: Los enunciados atómicos cuestionados que son incompatibles con su cuestionamiento son rechazados.

$$M(\neg) = M(I) = 1 \Rightarrow M(a\neg) = 0$$

Esta regla se ilustra en el paso 11 de la fig 21.

R¬I. *Rechazo del Cuestionamiento en la Incompatibilidad*: Si una fórmula atómica no es cuestionada entonces es incompatible con su cuestionamiento.

$$M(\neg) = 0 \Rightarrow M(I) = 1$$

Ra¬I. *Rechazo del alcance del Cuestionamiento en la Incompatibilidad*: Si una fórmula atómica es rechazada entonces es incompatible con su cuestionamiento.

$$M(a\neg) = 0 \Rightarrow M(I) = 1$$

Aa¬I. *Aceptación del alcance del Cuestionamiento en la Incompatibilidad*: Si una fórmula atómica es aceptada entonces es incompatible con su cuestionamiento.

$$M(a\neg) = 1 \Rightarrow M(I) = 1$$

nRa¬AI. *No rechazo del alcance del Cuestionamiento y Aceptación de la Incompatibilidad*: Los enunciados atómicos que no son rechazados y que son incompatibles con su cuestionamiento son aceptados.

$$M(a\neg) \neq 0 \text{ y } M(I) = 1 \Rightarrow M(a\neg) = 1$$

Esta regla se ilustra en el paso 11 de la fig 24 y también en el paso 10 de la fig 25.

Prueba de: A¬AI: Sean $M(\neg) = M(I) = 1$. Como $M(\neg) = 1$, por A¬ resulta que $M(a\neg) \neq 0$, es decir $M(a\neg) = 0$ o $M(a\neg) = *$. Pero, si $M(a\neg) = *$, por A*a¬I se infiere $M(I) = 0$, lo cual no es el caso. Por lo tanto, $M(a\neg) = 0$.

Prueba de R¬I: Sea $M(\neg) = 0$. Si $M(I) = 0$, por RI se deriva que $M(a\neg) = *$ y por nRa¬ resulta $M(\neg) = 1$, lo cual no es el caso.

Prueba de Ra¬I: Sea $M(a\neg) = 0$. Si $M(I) = 0$, por RI se deriva que $M(a\neg) = *$, lo cual no es el caso.

Prueba de Aa¬I: Sea $M(a\neg) = 1$. Si $M(I) = 0$, por RI se deriva que $M(a\neg) = *$, lo cual no es el caso.

Prueba de nRa¬AI: Sean $M(a\neg) \neq 0$ y $M(I) = 1$. Si $M(a\neg) = *$, por A*a¬I se infiere $M(I) = 0$, lo cual no es el caso, luego $M(a\neg) \neq *$, quedando sólo la alternativa $M(a\neg) = 1$.

## 5.3. Proposición. *Reglas derivadas para el condicional |*

nRiA→. *No Rechazo a la izquierda y Aceptación del Condicional*: Si el antecedente no es rechazado y el condicional es aceptado entonces el consecuente no es rechazado.

$$M(i\rightarrow) \neq 0 \text{ y } M(\rightarrow) = 1 \Rightarrow M(d\rightarrow) \neq 0$$

Esta regla se ilustra en el paso 10 de la fig 22. También en el paso 12.

RdA→. *Rechazo a la derecha y Aceptación del Condicional*: Si el consecuente es rechazado y el condicional es aceptado entonces el antecedente es rechazado.

$$M(d\rightarrow) = 0 \text{ y } M(\rightarrow) = 1 \Rightarrow M(i\rightarrow) = 0$$

Ri→. *Rechazo a la izquierda en el Condicional*: Si el antecedente es rechazado entonces el condicional es aceptado.

$$M(i\rightarrow) = 0 \Rightarrow M(\rightarrow) = 1$$

nRd→. *No Rechazo a la derecha en el Condicional*: Si el consecuente no es rechazado entonces el condicional es aceptado.

$$M(d\rightarrow) \neq 0 \Rightarrow M(\rightarrow) = 1$$

Prueba de nRiA→: Sean $M(i\rightarrow) \neq 0$ y $M(\rightarrow) = 1$. Supóngase que $M(d\rightarrow) = 0$, entonces por nRiRd→ se infiere $M(\rightarrow) = 0$, lo cual no es el caso. Por lo tanto, forzosamente $M(d\rightarrow) \neq 0$.

Prueba de RdA→: Sean M(d→) = 0 y M(→) = 1. Supóngase que M(i→) ≠ 0, entonces por nRiRd→ se infiere M(→) = 0, lo cual no es el caso. Por lo tanto, forzosamente M(i→) = 0.

Prueba de Ri→: Sea M(i→) = 0. Supóngase que M(→) = 0, entonces por R→ se infiere M(i) ≠ 0, lo cual no es el caso. Por lo tanto, forzosamente M(→) = 1.

Prueba de nRd→: Sea M(d→) ≠ 0. Supóngase que M(→) = 0, entonces por R→ se infiere M(d→) = 0, lo cual no es el caso. Por lo tanto, forzosamente M(→) = 1.

## 5.4. Proposición. *Reglas derivadas para la conjunción* |

nRiR∧. *No Rechazo a la izquierda y Rechazo de la Conjunción*: Si el coyunto izquierdo no es rechazado y la conjunción es rechazada entonces el coyunto derecho es rechazado.

$$M(i\wedge) \neq 0 \text{ y } M(\wedge) = 0 \Rightarrow M(d\wedge) = 0$$

nRdR∧. *No Rechazo a la derecha y Rechazo de la Conjunción*: Si el coyunto derecho no es rechazado y la conjunción es rechazada entonces el coyunto izquierdo es rechazado.

$$M(d\wedge) \neq 0 \text{ y } M(\wedge) = 0 \Rightarrow M(i\wedge) = 0$$

Rd∧. *Rechazo a la derecha en la Conjunción*: Si el coyunto derecho es rechazado entonces la conjunción es rechazada.

$$M(d\wedge) = 0 \Rightarrow M(\wedge) = 0$$

Ri∧. *Rechazo a la izquierda en la Conjunción*: Si el coyunto izquierdo es rechazado entonces la conjunción es rechazada.

$$M(i\wedge) = 0 \Rightarrow M(\wedge) = 0$$

Prueba de nRiR∧: Sean M(i∧) ≠ 0 y M(∧) = 0. Supóngase que M(d∧) ≠ 0, entonces por nRinRd∧ se infiere M(∧) = 1, lo cual no es el caso. Por lo tanto, forzosamente M(d∧) = 0.

Prueba de nRdR∧: Sean M(d∧) ≠ 0 y M(∧) = 0. Supóngase que M(i∧) ≠ 0, entonces por nRinRd∧ se infiere M(∧) = 1, lo cual no es el caso. Por lo tanto, forzosamente M(i∧) = 0.

Prueba de Rd∧: Sea M(d∧) = 0. Supóngase que M(∧) = 1, entonces por A∧ se infiere M(d∧) ≠ 0, lo cual no es el caso. Por lo tanto, forzosamente M(∧) = 0.

Prueba de Ri∧: Sea M(i∧) = 0. Supóngase que M(∧) = 1, entonces por A∧ se infiere M(i∧) ≠ 0, lo cual no es el caso. Por lo tanto, forzosamente M(∧) = 0.

## 5.5. Proposición. *Reglas derivadas para la disyunción* |

RdA∨. *Rechazo a la derecha y Aceptación de la Disyunción*: Si el disyunto derecho es rechazado y la disyunción es aceptada entonces el disyunto izquierdo no es rechazado.

$$M(d\vee) = 0 \text{ y } M(\vee) = 1 \Rightarrow M(i\vee) \neq 0$$

RiA∨. *Rechazo a la izquierda y Aceptación de la Disyunción*: Si el disyunto izquierdo es rechazado y la disyunción es aceptada entonces el disyunto derecho no es rechazado.

$$M(i\vee) = 0 \text{ y } M(\vee) = 1 \Rightarrow M(d\vee) \neq 0$$

Esta regla se ilustra en el paso 14 fig 24.

nRi∨. No *Rechazo a la izquierda en la Disyunción*: Si el disyunto izquierdo no es rechazado entonces la disyunción es aceptada.

$$M(i\vee) \neq 0 \Rightarrow M(\vee) = 1$$

nRd∨. No *Rechazo a la derecha en la Disyunción*: Si el disyunto derecho no es rechazado entonces la disyunción es aceptada.

$$M(d\vee) \neq 0 \Rightarrow M(\vee) = 1$$

Prueba de RdA∨: Sean M(d∨) = 0 y M(∨) = 1. Supóngase que M(i∨) = 0, entonces por RiRd∨ se infiere M(∨) = 0, lo cual no es el caso. Por lo tanto, forzosamente M(i∨) ≠ 0.

Prueba de RiA∨: Sean M(i∨) = 0 y M(∨) = 1. Supóngase que M(d∨) = 0, entonces por RiRd∨ se infiere M(∨) = 0, lo cual no es el caso. Por lo tanto, forzosamente M(d∨) ≠ 0.

Prueba de nRi∨: Sea M(i∨) ≠ 0. Supóngase que M(∨) = 0, entonces por R∨ se infiere M(i∨) = 0, lo cual no es el caso. Por lo tanto, forzosamente M(∨) = 1.

Prueba de nRd∨: Sea M(d∨) ≠ 0. Supóngase que M(∨) = 0, entonces por R∨ se infiere M(d∨) = 0, lo cual no es el caso. Por lo tanto, forzosamente M(∨) = 1.

## 5.6. Proposición. *Reglas derivadas para la negación fuerte* |

R~. *Rechazo de la Negación fuerte*: Una fórmula no es rechazada cuando se rechaza su negación.

$$M(\sim) = 0 \Rightarrow M(a\sim) \neq 0$$

nRa~. No *Rechazo del alcance de la Negación fuerte*: Una negación es rechazada cuando no se rechaza su alcance.

$$M(a\sim) \neq 0 \Rightarrow M(\sim) = 0$$

Prueba de R~: Sea M(~) = 0. Supóngase que M(a~) = 0, entonces por Ra~ se infiere M(~) = 1, lo cual no es el caso. Por lo tanto, forzosamente M(a~) ≠ 0.

Prueba de nRa~: Sea M(a~) ≠ 0. Supóngase que M(~) = 1, entonces por A~ se infiere M(a~) = 0, lo cual no es el caso. Por lo tanto, forzosamente M(~) = 0.

## **5.7. Proposición.** *Reglas derivadas para las marcas de atómicas* |

nAnRat. *No Aceptación y no Rechazo de atómicas*: Si una fórmula atómica no es aceptada y no es rechazada entonces es aceptada por defecto.

$$M(p) \neq 1 \text{ y } M(p) \neq 0 \Rightarrow M(p) = *$$

Esta regla se ilustra en el paso 9 de la fig 21.

nA*nRat. *No Aceptación por defecto y no Rechazo de atómicas*: Si una fórmula atómica no es aceptada por defecto y no es rechazada entonces es aceptada.

$$M(p) \neq 0 \text{ y } M(p) \neq * \Rightarrow M(p) = 1$$

nAnA*at. *No Aceptación y no Aceptación por defecto de atómicas*: Si una fórmula atómica no es aceptada y no es aceptada por defecto entonces es rechazada.

$$M(p) \neq 1 \text{ y } M(p) \neq * \Rightarrow M(p) = 0$$

Pruebas: Basta tener en cuenta que, por definición de función de marca, si p es atómica entonces M(p) = 1 o M(p) = 0 o M(p) = *.

## **5.8. Proposición.** *Reglas derivadas para las marcas de no atómicas* |

nAnat. *No Aceptación de no atómicas*: Si una fórmula no atómica (compuesta) no es aceptada entonces es rechazada.

$$M(X) \neq 1 \Rightarrow M(X) = 0$$

nRnat. *No Rechazo de no atómicas*: Si una fórmula no atómica (compuesta) no es rechazada entonces es aceptada.

$$M(X) \neq 0 \Rightarrow M(X) = 1$$

**Convención**: Estas reglas utilizan de forma *implícita*. Esto se ilustra en los pasos 2 y 4 de la fig 21, en los pasos 2, 4, 6 y 12 de la fig 22 y, en los pasos 2, 4, 9, 10 y 14 de la fig 24.

Pruebas: Basta tener en cuenta que, por definición de función de marca, si X es no atómica entonces M(X) = 1 o M(X) = 0.

## **5.9. Proposición.** *Regla de iteración* |

IM. *Iteración de Marca*: Sean n y k dos nodos asociados a una misma fórmula, si el nodo n está marcado entonces el nodo k tiene la misma marca.

[n asociado a W, k asociado a W, n ≠ k y M(n) conocida] ⇒ M(k) = M(n).

Esta regla se ilustra en el paso 10 de la fig 21. También en el paso 9 de la fig 22.

Prueba: Si n y k son nodos asociados a una misma fórmula W entonces n y k son ambas raíces de Ar[W], como m es una función entonces las marcas de n y k no pueden ser diferentes, es decir, M(n) = M(k).

## **5.10. Proposición.** *Reglas de opciones para la doble marca* |

OI-DM. *Opción Inicial que genera Doble Marca*: Si al suponer que un nodo n tiene una marca inicial y al aplicar las reglas para marcar nodos, se tiene como consecuencia marcas diferentes en algún par de nodos asociados a una misma fórmula, entonces el nodo n realmente no tiene la marca inicial.

Para cada nodo n, [M(n) = i ⇒ para algún nodo k, M(k) ≠ M(k)] ⇒ M(n) ≠ i.

RR-DM. *Rechazo de la Raíz que genera Doble Marca*: Si en el árbol de la fórmula X se supone que la raíz está marcada con 0 y al aplicar las reglas para marcar nodos, se tiene como consecuencia marcas diferentes en algún par de nodos asociados a una misma fórmula, entonces la fórmula X es A-válida. En este caso se dice que el árbol es un *árbol mal marcado*.

Para m una función de marca, [M(R[X]) = 0 ⇒ para algún nodo k, M(k) ≠ M(k)] ⇒ X es A-válida.

Observar que el Rechazo de la raíz (RR) es realmente una opción de rechazo.

Esta regla se ilustra en el paso 12 fig 21.

Prueba de OI-DM: Se tiene que M(n) = i ⇒ para algún nodo k, M(k) ≠ M(k). Supóngase que M(n) = i, se infiere entonces que para algún nodo k, M(k) ≠ M(k), pero esto es imposible ya que M es una función. Por lo tanto, forzosamente M(n) ≠ i.

Prueba de ORR-DM: Se tiene que M(R[X]) = 0 ⇒ para algún nodo k, M(k) ≠ M(k), donde m es una función de marca. Supóngase que X no es A-válido, es decir, existe una función de marca m, tal que M(R[X]) = 0, se infiere entonces que para algún nodo k, M(k) ≠ M(k), lo cual es imposible ya que M es una función. Por lo tanto, forzosamente X es A-válido.

## 5.11. Proposición. *Reglas de opciones para el condicional* |

OnRi-nRd→. *Opción de no Rechazo a la Izquierda que genera no Rechazo a la Derecha en un Condicional*: Si se supone que el antecedente no es rechazado y al aplicar las reglas para marcar nodos, se tiene como consecuencia que el consecuente no es rechazado, entonces el condicional realmente es aceptado.

$$[M(i\rightarrow) \neq 0 \Rightarrow M(d\rightarrow) \neq 0] \Rightarrow M(\rightarrow) = 1$$

ORd-Ri→. *Opción de Rechazo a la Derecha que genera Rechazo a la Izquierda en un Condicional*: Si se supone que el consecuente es rechazado y al aplicar las reglas para marcar nodos, se tiene como consecuencia que el antecedente es rechazado, entonces el condicional realmente es aceptado.

$$[M(d\rightarrow) = 0 \Rightarrow M(i\rightarrow) = 0] \Rightarrow M(\rightarrow) = 1$$

Esta regla se ilustra en el paso 10 fig 23.

Prueba de OnRi-nRd→: Se tiene que M(i→) ≠ 0 ⇒ M(d→) ≠ 0. Supóngase que M(→) ≠ 1, es decir M(→) = 0, entonces por R→ resulta que M(i→) ≠ 0 y M(d→) = 0, pero al tener M(i→) ≠ 0 se infiere que M(d→) ≠ 0, pero esto es imposible ya que se tiene que M(d→) = 0. Por lo tanto, forzosamente M(→) = 1.

Prueba de ORd-Ri→: Se tiene que M(d→) = 0 ⇒ M(i→) = 0. Supóngase que M(→) ≠ 1, es decir M(→) = 0, entonces por R→ resulta que M(i→) ≠ 0 y M(d→) = 0, pero al tener M(d→) = 0 se infiere que M(i→) = 0, pero esto es imposible ya que se tiene que M(i→) ≠ 0. Por lo tanto, forzosamente M(→) = 1.

## 5.12. Proposición. *Reglas de opciones para la disyunción* |

ORi-nRd∨. *Opción de Rechazo a la Izquierda que genera no Rechazo a la Derecha en una Disyunción*: Si se supone que el disyunto izquierdo es rechazado y al aplicar las reglas para marcar nodos, se tiene como consecuencia que el disyunto derecho no es rechazado, entonces la disyunción realmente es aceptada.

$$[M(i\vee) = 0 \Rightarrow M(d\vee) \neq 0] \Rightarrow M(\vee) = 1$$

Esta regla se ilustra en el paso 4 de la fig 26.

ORd-nRi∨. *Opción de Rechazo a la Derecha que genera no Rechazo a la Izquierda en una Disyunción*: Si se supone que el disyunto derecho es rechazado y al aplicar las reglas para marcar nodos, se tiene como consecuencia que el disyunto izquierdo no es rechazado, entonces la disyunción realmente es aceptada.

$$[M(d\vee) = 0 \Rightarrow M(i\vee) \neq 0] \Rightarrow M(\vee) = 1$$

Prueba de ORi-nRd∨: Se tiene que M(i∨) = 0 ⇒ M(d∨) ≠ 0. Supóngase que M(∨) ≠ 1, es decir M(∨) = 0, entonces por R∨ resulta que M(i∨) = 0 y M(d∨) = 0, pero al tener M(i∨) = 0 se infiere que M(d∨) ≠ 0, pero esto es imposible ya que se tiene que M(d∨) = 0. Por lo tanto, forzosamente M(∨) = 1.

Prueba de ORd-nRi∨: Se tiene que M(d∨) = 0 ⇒ M(i∨) ≠ 0. Supóngase que M(∨) ≠ 1, es decir M(∨) = 0, entonces por R∨ resulta que M(i∨) = 0 y M(d∨) = 0, pero al tener M(d∨) = 0 se infiere que M(i∨) ≠ 0, pero esto es imposible ya que se tiene que M(i∨) = 0. Por lo tanto, forzosamente M(∨) = 1.

## 5.13. Proposición. *Unicidad de la extensión de una función de marca* |

Cada función de marca de hojas *m*, puede ser extendida de manera única, a una *función de marca de nodos, M*, de N(X) en {0, 1, *}, haciendo *M*(h) = *m*(h) si h es una hoja, y aplicando las reglas de instanciación de espacios vacíos, junto con las reglas primitivas y derivadas para el forzamiento de marca, las cuales son presentadas en las secciones 4.1 a 4.7 y 5.1 a 5.12.

Prueba: Se debe probar que, a) la extensión M de *m* se aplique a todas las fórmulas, b) la asignación de M a cada fórmula sea única, y c) no existe otra extensión M' de *m*, la cual se aplica a todas las fórmulas.

Parte a. Por la definición, la extensión M de *m* se aplica a todas las fórmulas.

Parte b. Se prueba por inducción sobre la profundidad P del árbol de las fórmulas.

Paso base. Profundidad 0, significa que el nodo es una hoja, en este caso M coincide con la función *m*, por lo que se satisface la unicidad de asignación.

Paso inductivo. Sea P(Ar[X]=L, con L>0. Como hipótesis inductiva se tienen:

P(Ar[Y]<L implica la asignación de M a Y es única, P(Ar[Z]<L implica la asignación de M a Z es única, P(Ar[W*a*]<L implica la asignación de M a W*a* es única.

Caso 1. X=~Z. Supóngase que M(~)=1 y M(~)=0. Como M(~)=1, por la regla A~ se deriva que M(a~)=0, además, como, por la regla R~ se deriva que M(a~)=1, resultando que la asignación de M a Ar[Z] no es única, lo cual contradice la hipótesis inductiva. Por lo tanto, la asignación de M a ~Z es única.

Caso 2. X=Y∧Z. Supóngase que M(∧)=1 y que M(∧)=0. Como M(∧)=1, por la regla A∧ se obtiene que M(i∧)≠0 y M(d∧)≠0, se tiene entonces que M(∧)=0 y M(i∧)≠0, aplicando la regla RiA∧ se infiere que M(d∧)=0, pero M(d∧)≠0, lo cual contradice la hipótesis inductiva. Por lo tanto, la asignación de M a Y∧Z es única.

Caso 3. X=Y∨Z. Supóngase que M(∨)=1 y que M(∨)=0. Como M(∨)=0, por la regla R∨ se deriva que M(i∨)=0 y M(d∨)=0, se tiene entonces que M(∨)=1 y M(i∨)=0, aplicando la regla RiA∨ se infiere que M(d∨)≠0, pero M(d∨)=0, lo cual contradice la hipótesis inductiva. Por lo tanto, la asignación de M a Y∨Z es única.

Caso 4. X=Y→Z. Supóngase que M(→)=1 y que M(→)=0. Como M(→)=0, por la regla R→ se deriva que M(i→)≠0 y M(d→)=0, se tiene entonces que M(→)=1 y M(i→)≠0, aplicando la regla nRiA→ se infiere que M(d→)≠0, pero M(d→)=0, lo cual contradice la hipótesis inductiva. Por lo tanto, la asignación de M a Y→Z es única.

Caso 5. X=¬Z. Supóngase que M(¬)=1 y M(¬)=0. Como M(¬)=0, por la regla A~ se deriva que M(a¬)≠1, además, como M(¬)=0, por la regla R¬ se deriva que M(a¬)=1, resultando que la asignación de M a Ar[a¬] no es única, lo cual contradice la hipótesis inductiva. Por lo tanto, la asignación de M a ¬Z es única.

Caso 6. X=IZ. Supóngase que M(I)=1 y M(I)=0. Como M(I)=0, por la regla RI se deriva que M(a¬)=*, es decir M(a¬)≠0, además, como M(I)=0, por la regla nRa¬AI se deriva que M(a¬)=1, resultando que la asignación de M a Ar[I] no es única, lo cual contradice la hipótesis inductiva. Por lo tanto, la asignación de M a IZ es única.

Por el principio de inducción matemática, se ha probado que la asignación de M a cada fórmula es única.

Parte c. Supóngase que existe otra extensión M' de *m*, la cual es una función que se aplica a todas las fórmulas. Si M≠M', entonces existe al menos una fórmula F, tal que M(Ar[F])≠M'(Ar[F]). Entre estas fórmulas, el principio del buen orden, garantiza la existencia de al menos una fórmula de profundidad mínima, sea X una de estas fórmulas, por lo que P(Ar[X])=L y es mínima, es decir, para cada fórmula T, si P(Ar[T])<L entonces M(Ar[T])=M'(Ar[T]). Además, como M y M' son ambas extensiones de *m*, entonces Ar[X] no puede ser una hoja, por lo que X debe ser una fórmula compuesta.

Caso 1. X=~Z. Supóngase que M(~)=1 y que M'(~)=0. Como M(~)=1, aplicando la regla A~ se infiere que M(a~)=0, además, como M'(~)=0, utilizando la regla R~ se infiere que M'(a~)=1, por lo que M(a~)≠ M'(a~), lo cual no es posible puesto que P(a~)<L.

Caso 2. X=Y∧Z. Supóngase que M(∧)=1 y que M'(∧)=0. Como M(∧)=1, aplicando la regla A∧ se infiere que M(i∧)≠0 y M(d∧)≠0, además P(Ar[i∧])<L, de donde M'(i∧)≠1, se tiene entonces que M'(∧)=0 y M'(i∧)≠0, utilizando la regla nRiR∧ se infiere que M'(d∧])=0, por lo que M(d∧)≠ M'(d∧), lo cual no es posible puesto que P(d∧)<L.

Caso 3. X=Y∨Z. Supóngase que M(∨)=1 y que M'(∨)=0. Como M(∨])=0, aplicando la regla R∨ se infiere que M(i∨)=0 y M(d∨)=0, además P(d∨)<L, de donde M'(d∨)=0, se tiene entonces que M'(∨)=1 y M'(d∨)=0 utilizando la regla RdA∨ se infiere que M'(i∨)=1, por lo que M(i∨)≠ M'(i∨), lo cual no es posible puesto que P(i∨)<L.

Caso 4. X=Y→Z. Supóngase que M(→)=1 y que M'(→)=0. Como M(→)=0, aplicando la regla R→ se infiere que M(i→)≠0 y M(d→)=0, además P(i→)<L, de donde M'(i→)≠0, se tiene entonces que M'(i→)≠0 y M'(→)=1 utilizando la regla nRiA→ se infiere que M'(d→)≠0, por lo que M(d→)≠ M'(d→), lo cual no es posible puesto que P(d→)<L.

Caso 5. X=¬Z. Supóngase que M(¬)=1 y que M'(¬)=0. Como M(¬)=1, aplicando la regla A¬ se infiere que M(a¬)≠1, además, como M'(¬)=0, utilizando la regla R¬ se infiere que M'(a¬)=1, por lo que M(a¬)≠ M'(a¬), lo cual no es posible puesto que P(a¬)<L.

Caso 6. X=IZ. Supóngase que M(I)=1 y que M'(I)=0. Como M(I)=0, aplicando la regla RI se infiere M(a¬)=*, es decir M(a¬)≠0, además P(a¬)<L, de donde M'(a¬)≠0, se tiene entonces que M'(a¬)≠0 y M(I)=1 utilizando la regla nRa¬AI se infiere que M'(a¬)=1, por lo que M(a¬)≠ M'(a¬), lo cual no es posible puesto que P(a¬)<L.

Por lo tanto, no existe otra extensión M' de *m*, la cual se aplica a todas las fórmulas.

## 5.14. Presentación gráfica de las reglas para el forzamiento de marcas |

Un *nodo encerrado en un círculo* indica que el nodo está marcado con 1 (es aceptado), un *nodo encerrado en un cuadro* indica que el nodo está marcado con 0 (es rechazado), un *nodo encerrado en un triángulo* indica que el nodo está marcado con * (es aceptado por defecto), un *nodo encerrado en un doble cuadro* indica que el nodo no está marcado con 0 (no es rechazado), un *nodo encerrado en un círculo dentro de un cuadro* indica que el nodo no está marcado con 1 (no es aceptado), un *nodo encerrado en un triángulo dentro de un cuadro* indica que el nodo no está marcado con * (no es aceptado por defecto), un *nodo encerrado en una estrella* indica que el nodo está marcado con 0, 1 o *. Cuando las marcas están punteadas indican que es una opción (un supuesto), no es el resultado de aplicar una regla. La configuración inicial se muestra sobre la línea punteada y las marcas iniciales se presentan en color azul, la configuración que resulta cuando se aplica la regla, se presenta debajo de la línea punteada y la nueva marca generada por la regla se presenta en color rojo.

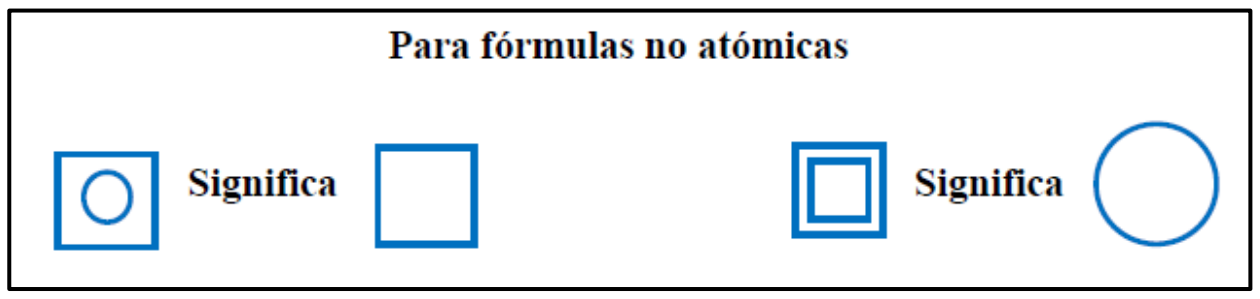

Figura 3. Reglas para las fórmulas no atómicas

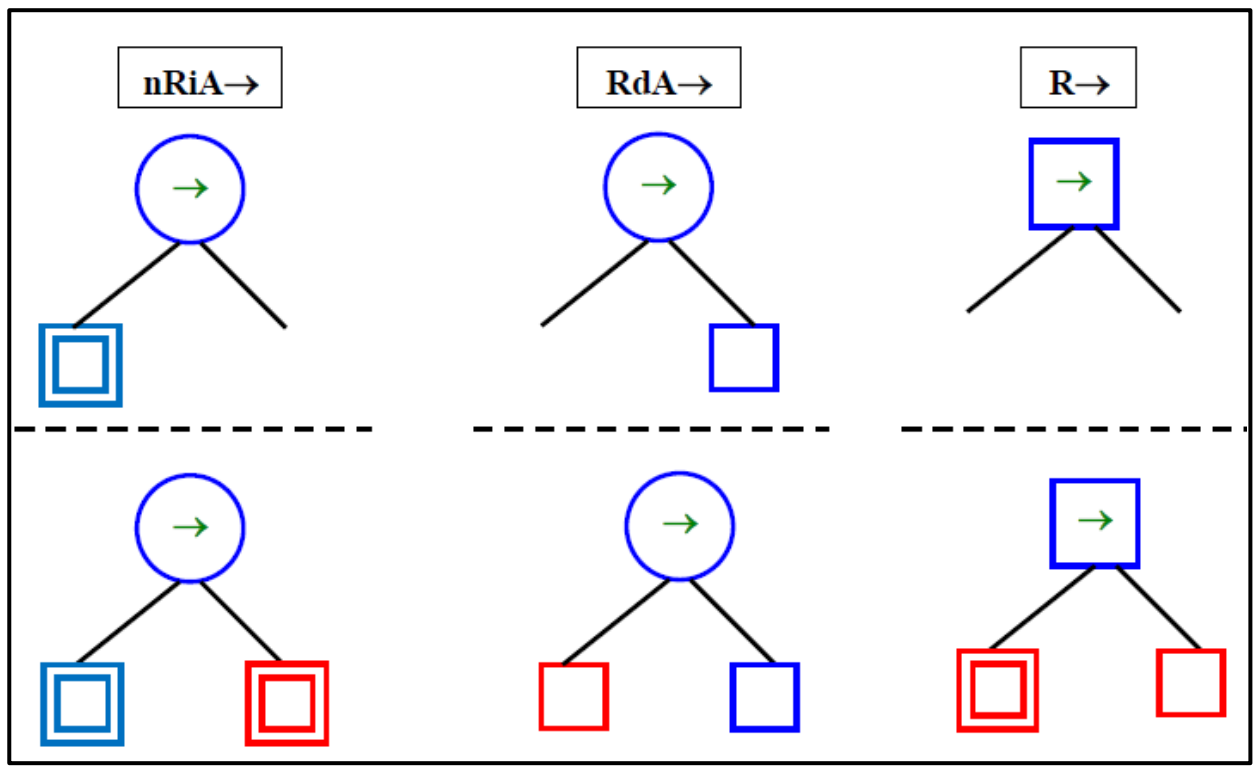

Figura 4. Reglas para el condicional

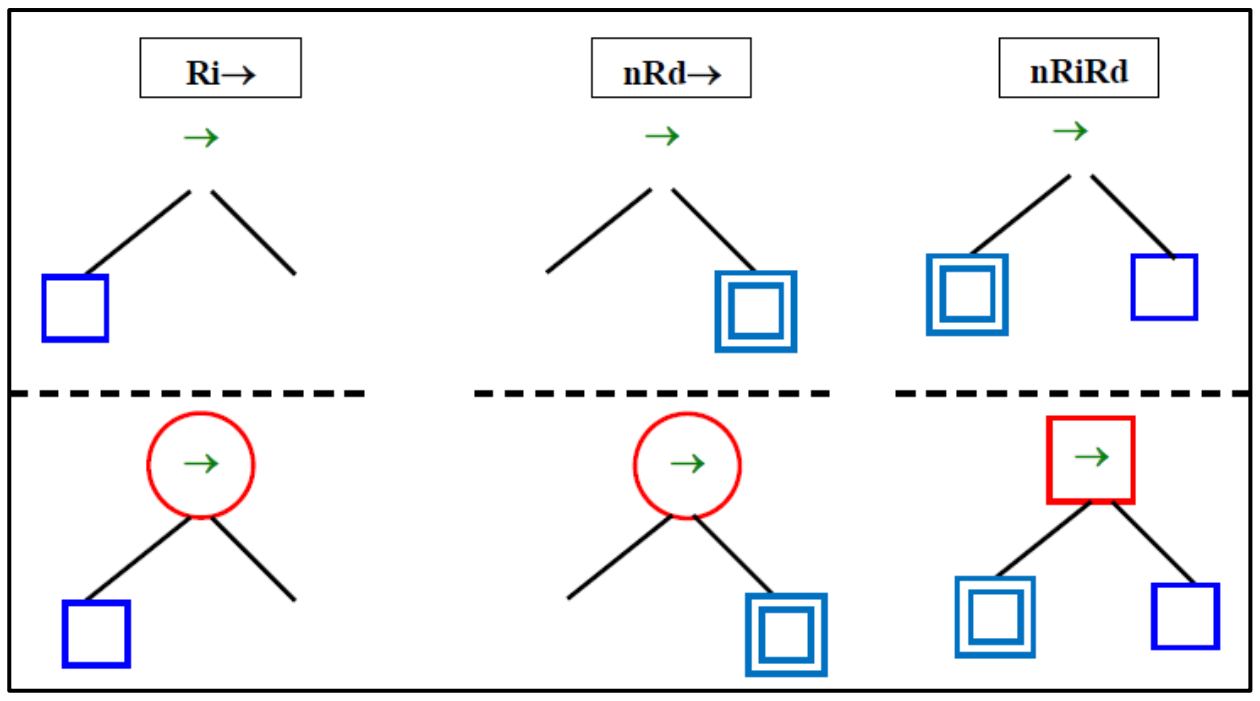

Figura 5. Reglas para el condicional

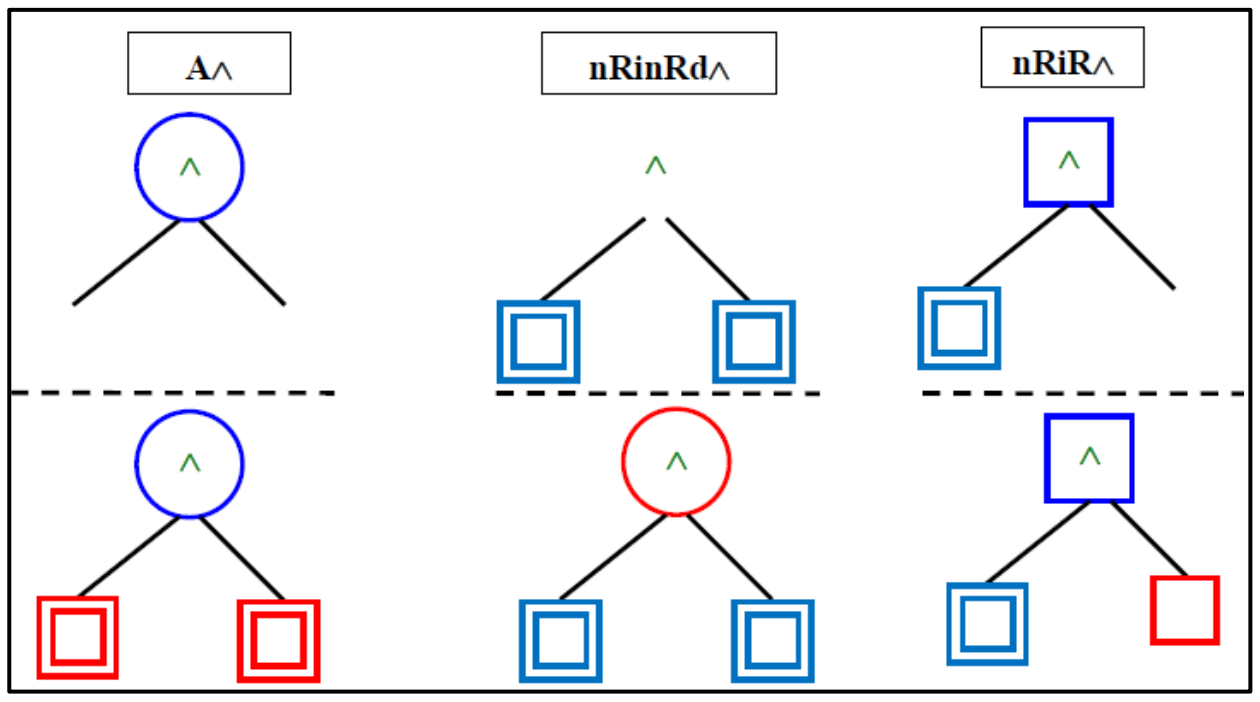

Figura 6. Reglas para la conjunción

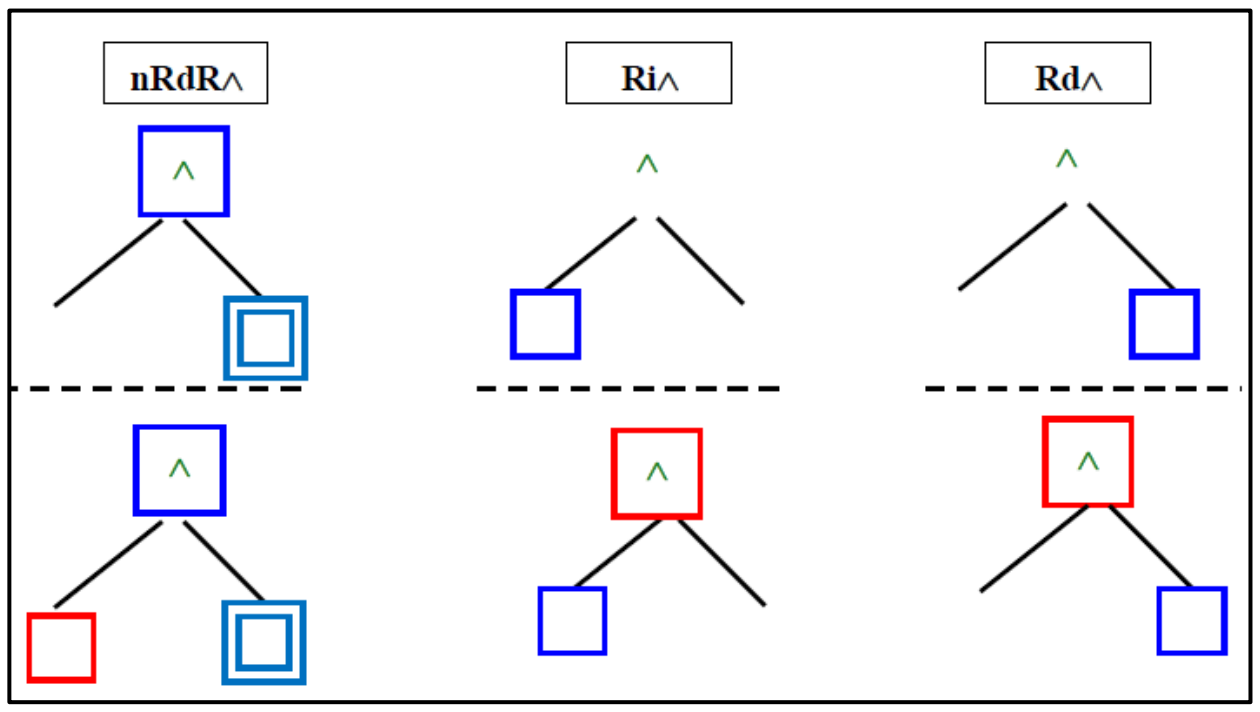

Figura 7. Reglas para la conjunción

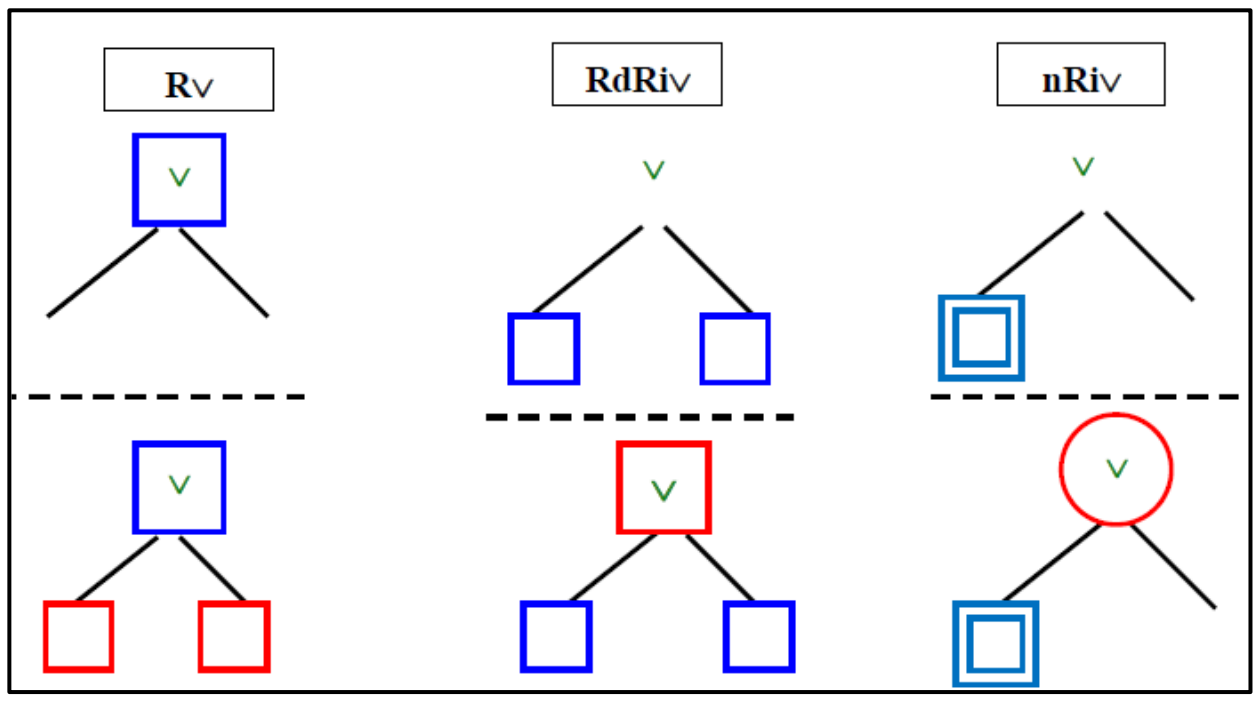

Figura 8. Reglas para la disyunción

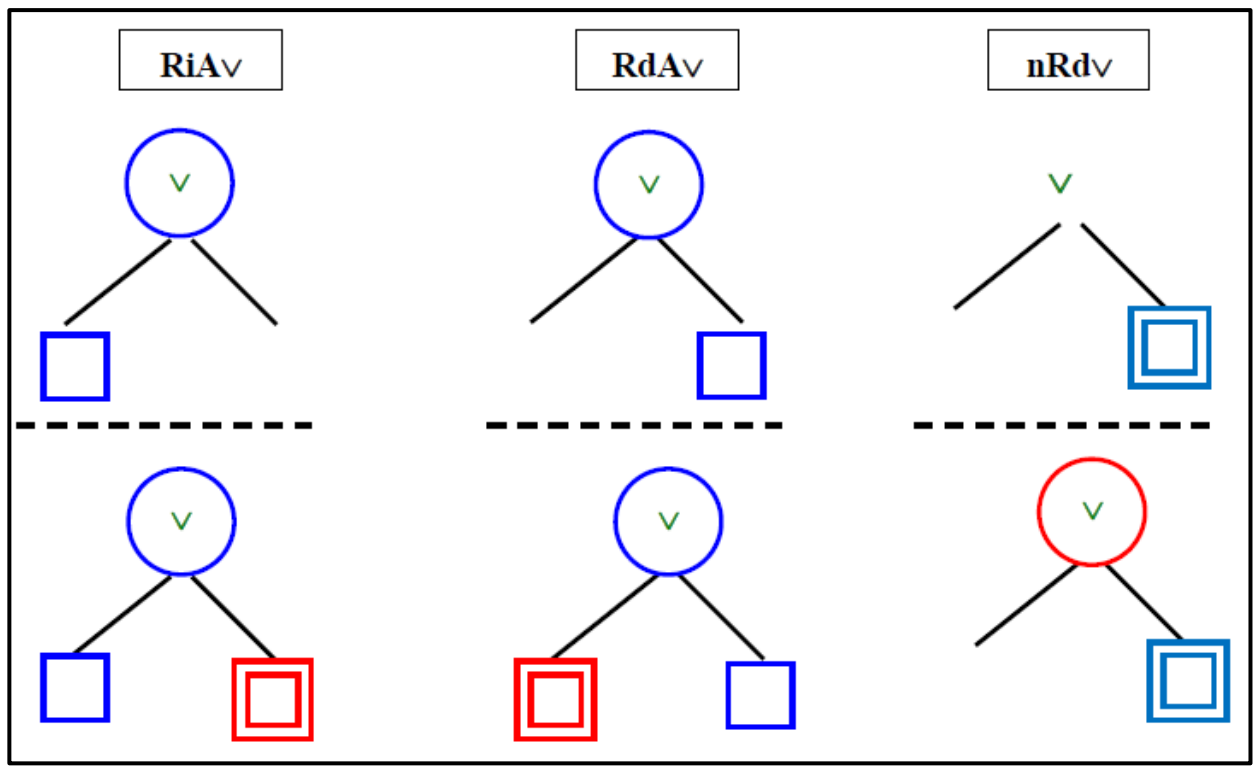

Figura 9. Reglas para la disyunción

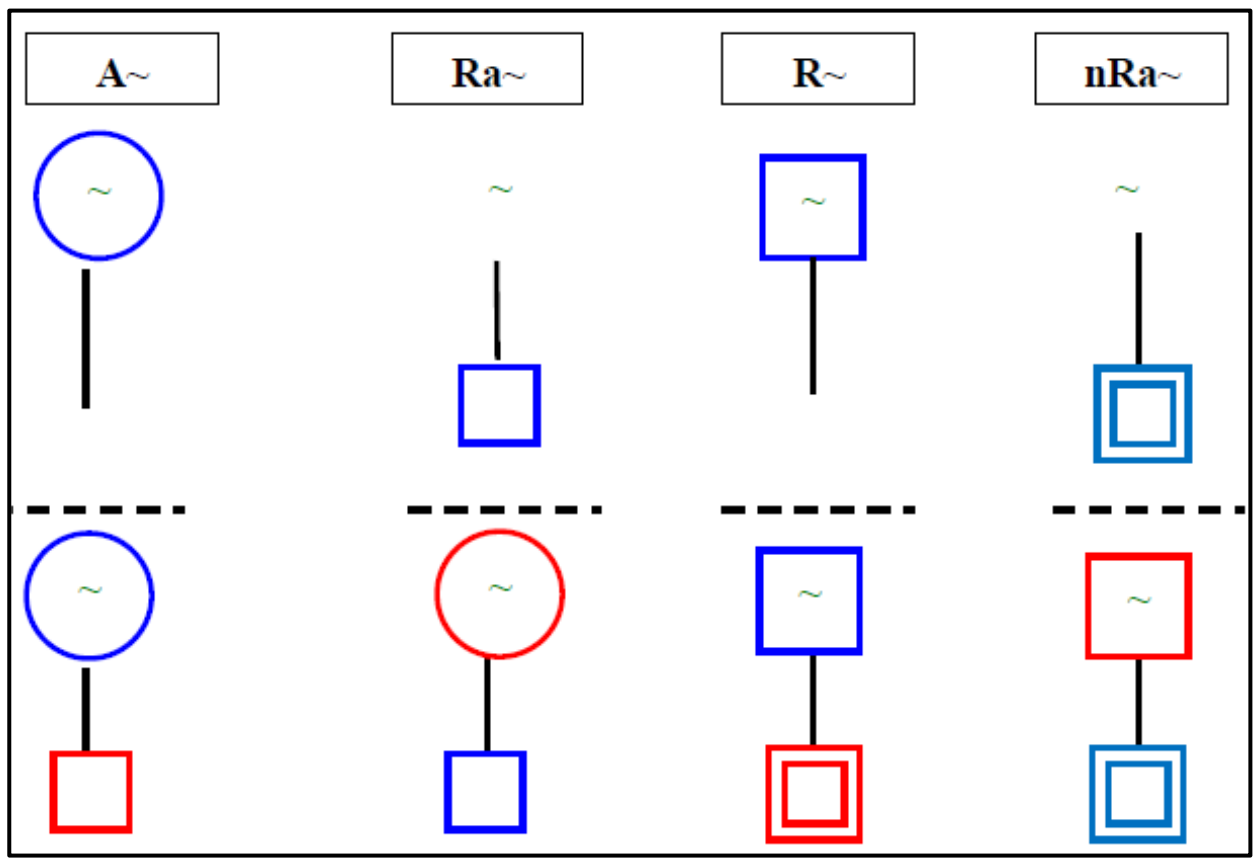

Figura 10. Reglas para la negación fuerte

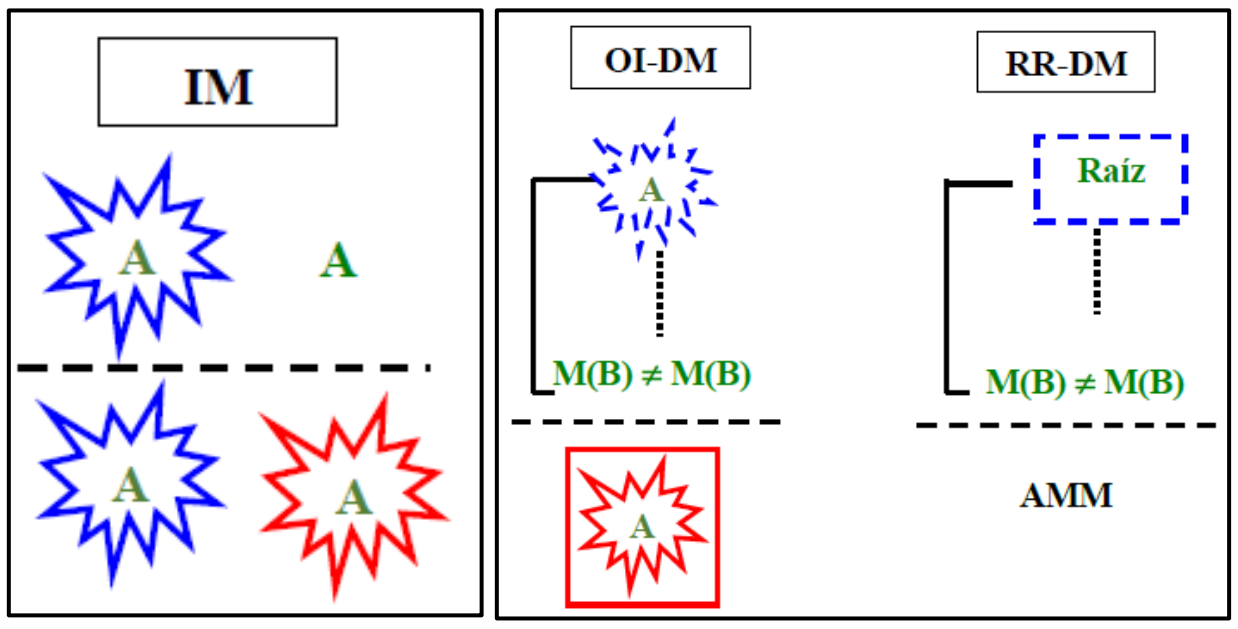

Figura 11. Reglas para la iteración y la doble marca

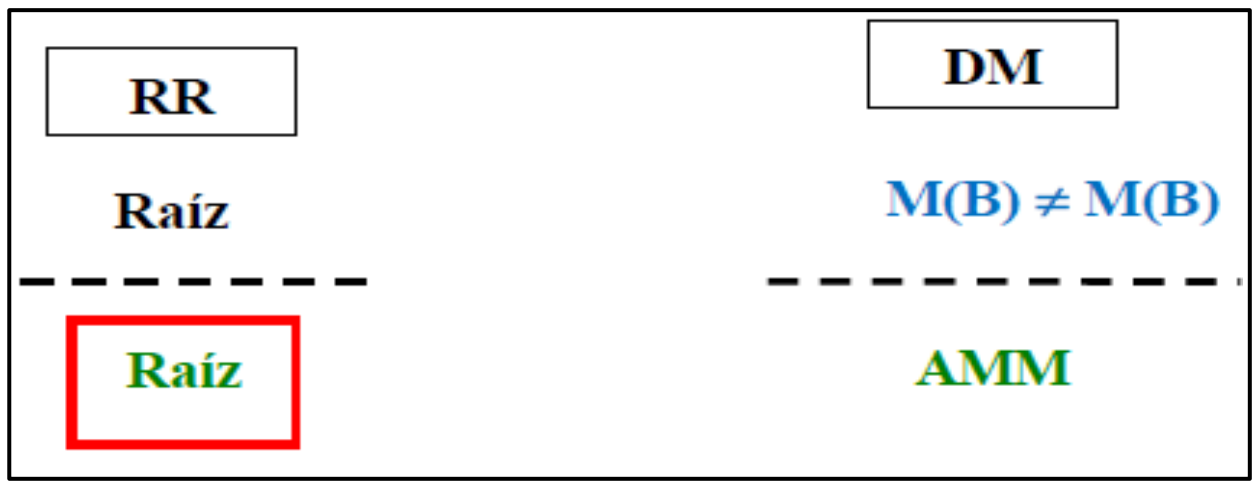

Figura 12. Reglas para el rechazo de la raíz

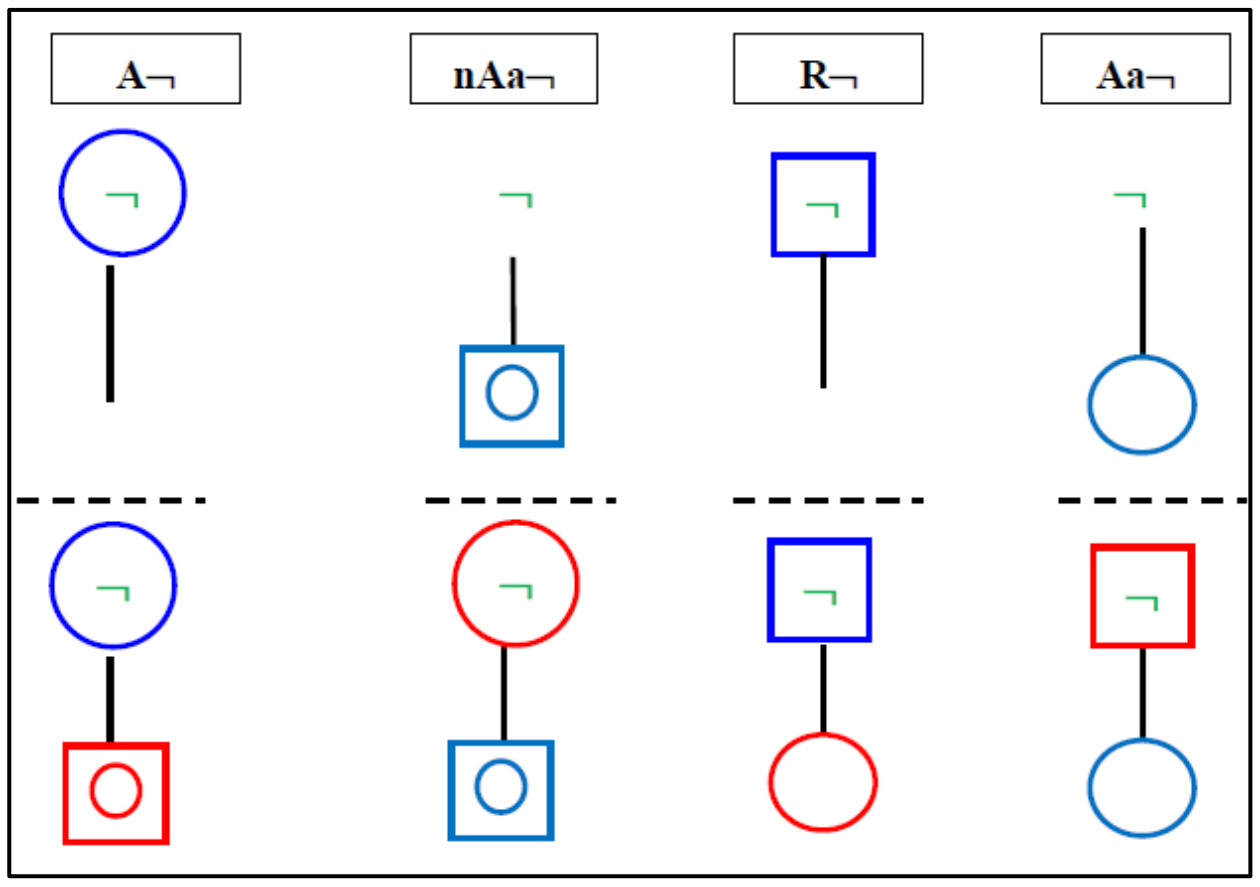

Figura 13. Reglas para la negación débil

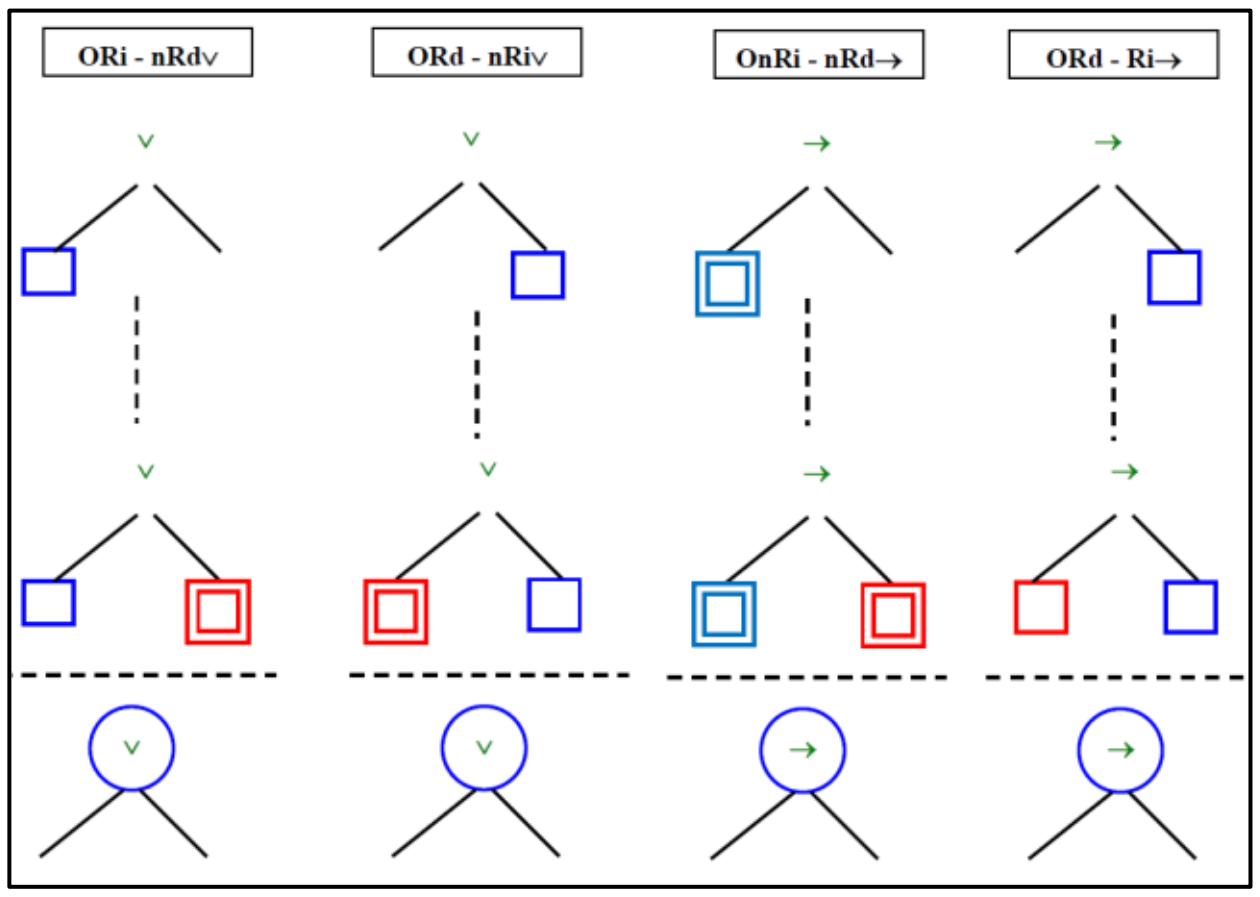

Figura 14. Reglas de opciones en disyunción y condicional

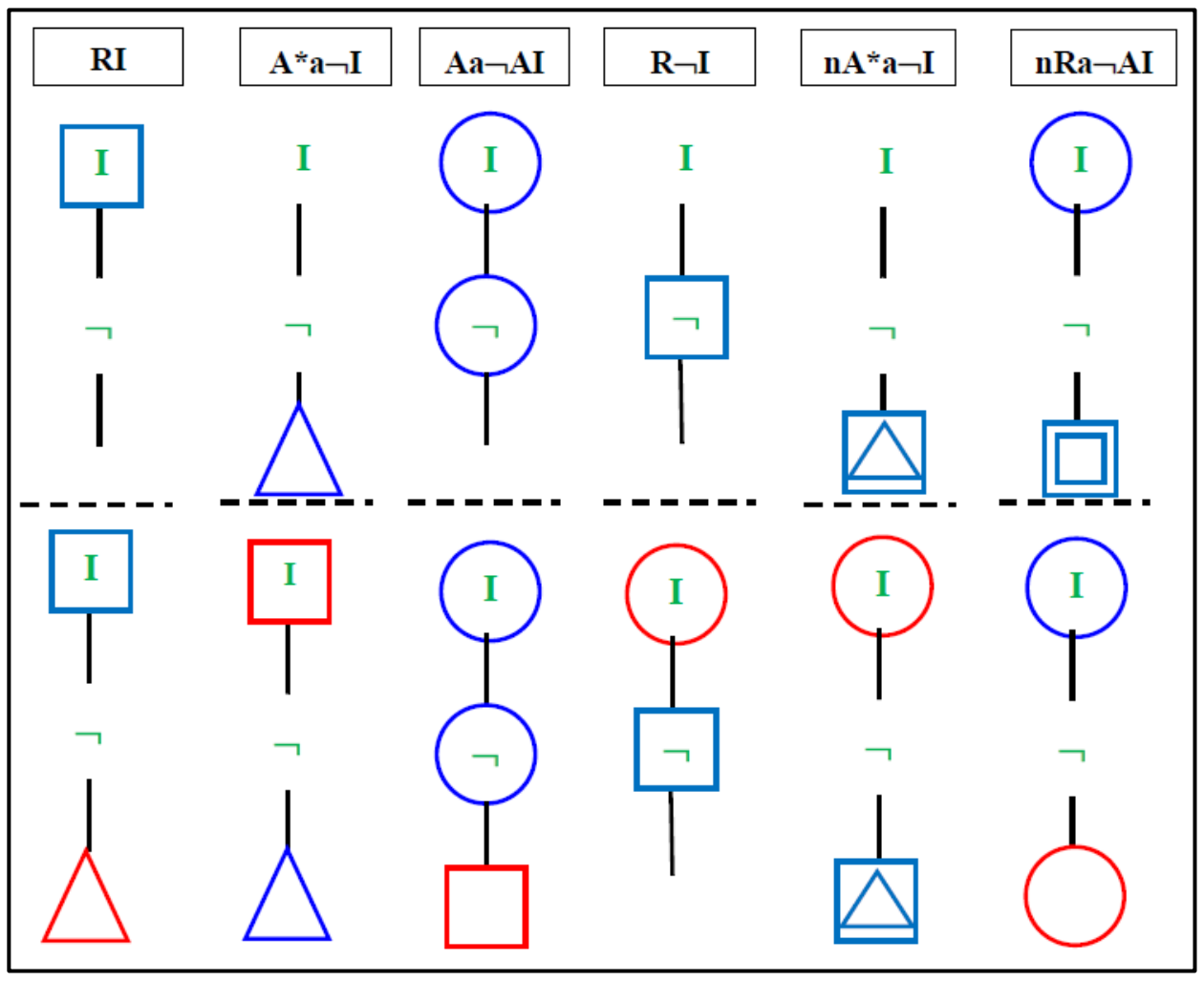

Figura 15. Reglas para la incompatibilidad

El bicondicional entre las fórmulas X y Y, denotado X↔Y, se define de la manera habitual como (X→Y)∧(Y→X). En consecuencia, las reglas para el bicondicional son las siguientes:

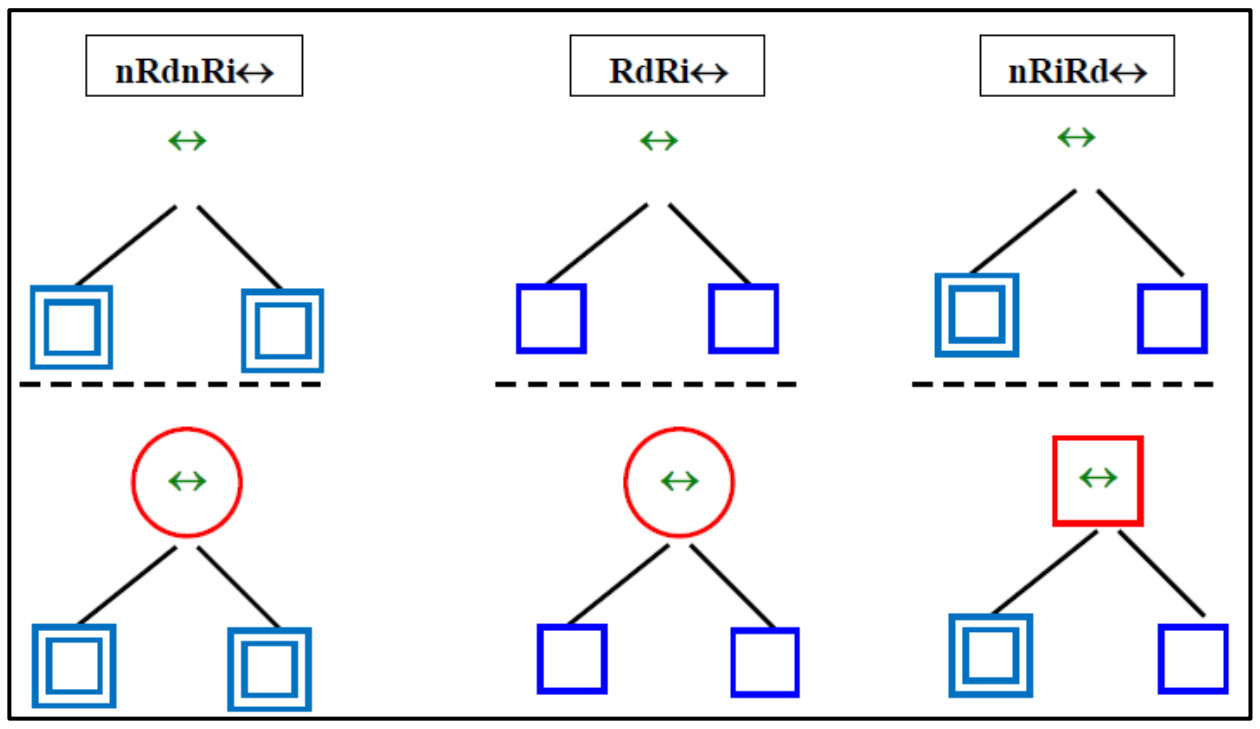

Figura 16. Reglas el bicondicional

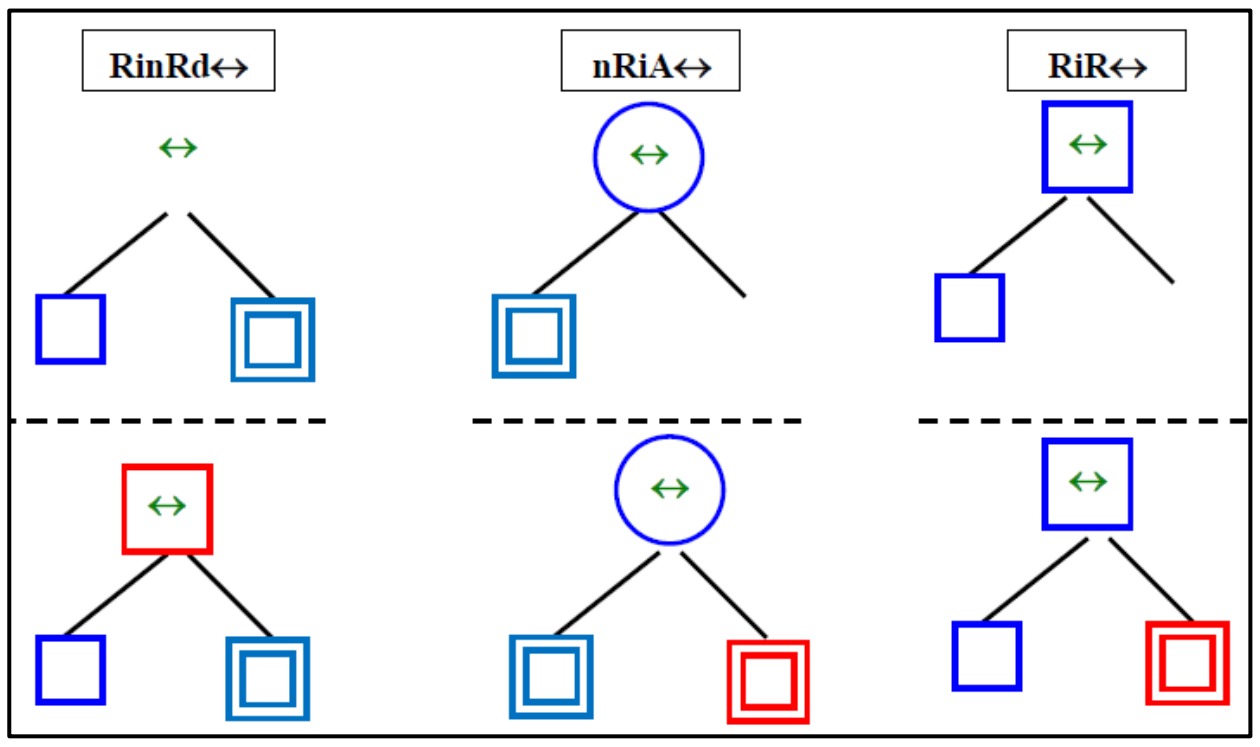

Figura 17. Reglas el bicondicional

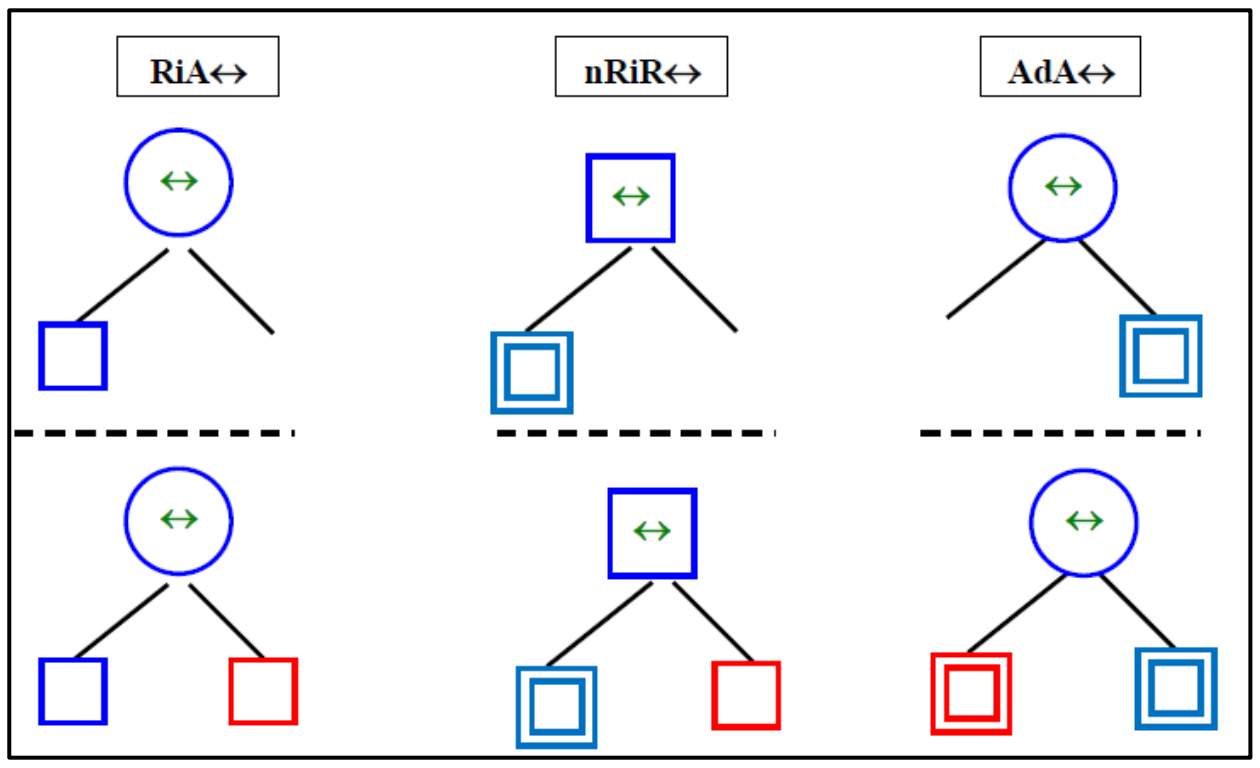

Figura 18. Reglas el bicondicional

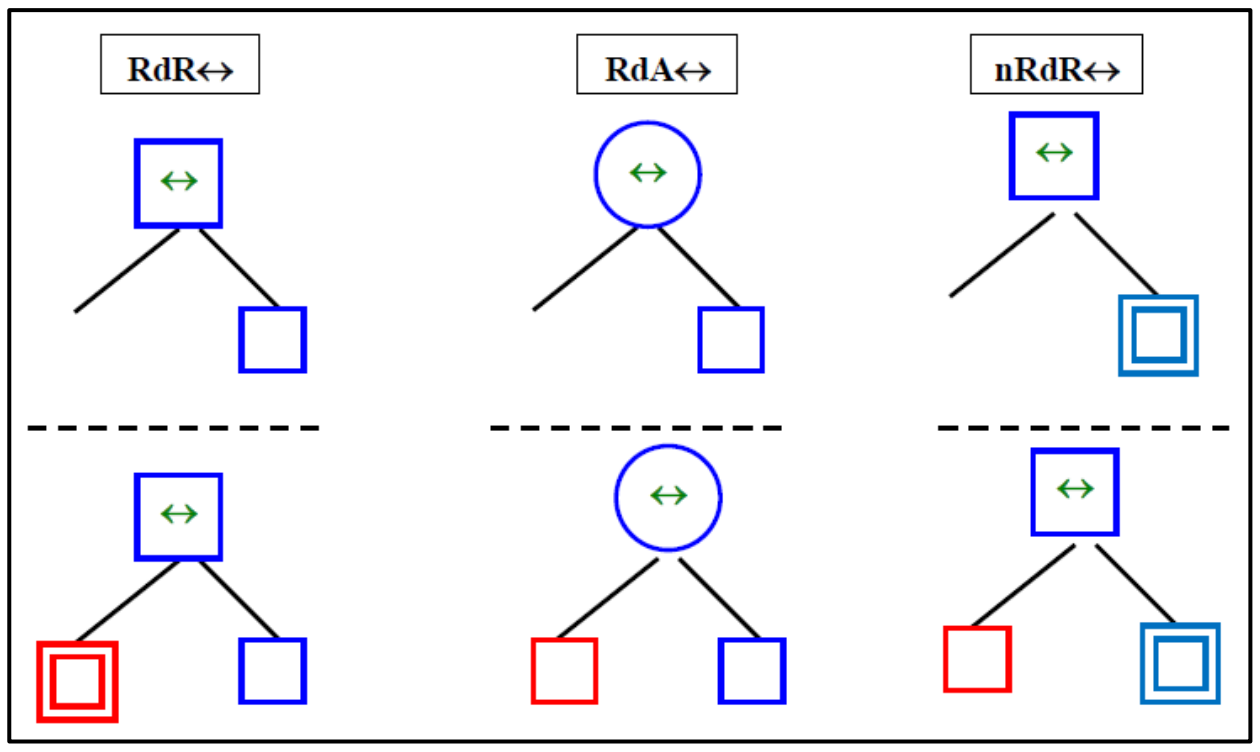

Figura 19. Reglas el bicondicional

## 6. SEMÁNTICA PARA EL SISTEMA DEDUCTIVO P1 |

En (Sette, 1973), se presenta una semántica tri-valuada para el sistema deductivo P1, mediante las siguientes tablas de verdad de la figura 20, donde 1 y * son los valores designados:

| → | 1 | * | 0 |
|---|---|---|---|
| 1 | 1 | 1 | 0 |
| * | 1 | 1 | 0 |
| 0 | 1 | 1 | 1 |

| ∧ | 1 | * | 0 |
|---|---|---|---|
| 1 | 1 | 1 | 0 |
| * | 1 | 1 | 0 |
| 0 | 0 | 0 | 0 |

| ∨ | 1 | * | 0 |
|---|---|---|---|
| 1 | 1 | 1 | 1 |
| * | 1 | 1 | 1 |
| 0 | 1 | 1 | 0 |

| | ¬ | ~ | I |
|---|---|---|---|
| 1 | 0 | 0 | 1 |
| * | 1 | 0 | 0 |
| 0 | 1 | 1 | 1 |

| ↔ | 1 | * | 0 |
|---|---|---|---|
| 1 | 1 | 1 | 0 |
| * | 1 | 1 | 0 |
| 0 | 0 | 0 | 1 |

Figura 20

Una *asignación*, v, es una función del conjunto de fórmulas atómicas en el conjunto {1, *, 0}.

La asignación v, se extiende a una función, V, del conjunto de fórmulas en el conjunto {1, *, 0}, de la siguiente manera:

V(p) = v(p) si p es atómica.

Para las fórmulas compuestas, la función V se encuentra definida por las tablas de verdad trivalentes de la figura 20.

Se dice que una fórmula X es T-*válida* (válida desde el punto de vista de la semántica trivalente), si y solamente si para toda asignación, v, V(X) ≠ 0.

## 7. EQUIVALENCIA ENTRE A-VALIDEZ Y T-VALIDEZ |

Se define la *Complejidad* C, como una función la cual asigna a cada fórmula de P1 un entero no negativo de la siguiente forma:

C(p) = 0, donde p es una fórmula atómica.

C(Ip) = 2, donde p es una fórmula atómica.

C(XkY) = 1 + máximo de {C(X), C(Y)}, donde k∈{∧, ∨, →}.

C(~X) = C(¬X) = 1 + C(X).

Se define la *Profundidad* P, como una función la cual asigna a cada nodo de un árbol un entero no negativo de la siguiente forma:

P(p) = 0, donde p es una fórmula atómica.

P(Ip) = 2, donde p es una fórmula atómica.

P(Ar[XkY]) = 1 + máximo de {P(Ar[X]), P(Ar[Y])}, donde k∈{∧, ∨, →}.

P(Ar[~X]) = P(Ar[¬X]) = 1 + P(Ar[X]).

## 7.1. Proposición. *Asignación correspondiente a una función de marca* |

Para cada fórmula X de P1 y para cada función de marca *m*, existe una asignación $v_m$, tal que, M(R[X]) = 0 ⇔ $V_m$(X) = 0.

Prueba: Sea p una fórmula atómica de P1, y sea *m* una función de marcas. Se define la función $v_m$ del conjunto de fórmulas atómicas en el conjunto {0, *, 1} de la siguiente forma:

$$v_m(p) = m(p)$$

La función $v_m$ se extiende a una función $V_m$ del conjunto de fórmulas de P1 en el conjunto {0, *, 1}, haciendo $V_m$(p) = $v_m$(p) cuando p es atómica, para las fórmulas compuestas $V_m$ se define mediante las *tablas de verdad trivalentes* de la sección 6. Se tiene entonces que $v_m$ es una asignación de P1.

Para probar que M(R[X]) = 0 ⇔ $V_m$(X) = 0, se procede por inducción sobre la *complejidad* de la fórmula X.

Paso base.

Supóngase que la C(X) = 0, esto significa que X es una fórmula atómica. Sea X = p, con p atómica. M(R[X]) = M(R[p]) = M(p) = m(p) = $v_m$(p) = $V_m$(p) = $V_m$(X), por lo tanto, M(R[X]) = 0 ⇔ $V_m$(X) = 0.

Paso de Inducción: Supóngase que C(X) ≥ 1, por lo que X tiene una de las siguientes formas: ¬p, ~p, Ip, ¬W, ~W, Y∧Z, Y∨Z o Y→Z, donde p es atómica y W no es atómica.

Caso 1: X = ¬p. M(R[¬p]) = 0 ⇔ M[¬] = 0, por las reglas R¬ y Aa¬ se tiene M[¬] = 0 ⇔ M[p] = 1, pero, por la definición de $v_m$, M[p] = m(p) = $v_m$(p) = $V_m$(p), y por la tabla de verdad se tiene $V_m$(p) = 1 ⇔ $V_m$(¬p) = 0. Resultando que M(R[¬p]) = 0 ⇔ $V_m$(¬p) = 0, es decir, M(R[X]) = 0 ⇔ $V_m$(X) = 0.

Caso 2: X = ~p. M(R[~p]) = 0 ⇔ M[~] = 0, por las reglas R~ y Aa~ se tiene M[~] = 0 ⇔ M[p] ≠ 0, pero, por la definición de $v_m$, M[p] = m(p) = $v_m$(p) = $V_m$(p), y por la tabla de verdad se tiene $V_m$(p) ≠ 0 ⇔ $V_m$(~p) = 0. Resultando que M(R[~p]) = 0 ⇔ $V_m$(~p) = 0, es decir, M(R[X]) = 0 ⇔ $V_m$(X) = 0.

Caso 3: X = Ip. M(R[Ip]) = 0 ⇔ M[I] = 0, por las reglas RI y A*a¬I se tiene M[I] = 0 ⇔ M[p] = *, pero, por la definición de $v_m$, M[p] = m(p) = $v_m$(p) = $V_m$(p), y por la tabla de verdad se tiene $V_m$(p) = * ⇔ $V_m$(Ip) = 0. Resultando que M(R[Ip]) = 0 ⇔ $V_m$(Ip) = 0, es decir, M(R[X]) = 0 ⇔ $V_m$(X) = 0.

Caso 4: Sea X = ¬W, con W no atómica. La hipótesis inductiva

M(R[W]) = 0 ⇔ $V_m$(W) = 0, al ser W no atómica equivale a M(R[W]) = 1 ⇔ $V_m$(W) = 1. Se tiene que R[X] = ¬, por lo que M(R[X]) = 0 ⇔ M(¬) = 0, pero por R¬ y Aa¬ se tiene M(¬) = 0 ⇔ M(a¬) = 1, y como además a¬ = R[W], resulta que M(R[X]) = 0 ⇔ [M(R[W]) = 1]. Utilizando la hipótesis inductiva se tiene que M(R[X]) = 0 ⇔ $V_m$(W) = 1. Por las tablas de verdad se concluye que M(R[X]) = 0 ⇔ $V_m$(¬W) = 0, es decir, M(R[X]) = 0 ⇔ $V_m$(X) = 0.

Caso 5: Sea X = ~W, con W no atómica. Se tiene que R[X] = ~, por lo que M(R[X]) = 0 ⇔ M(~) = 0, pero por R~ y Aa~ se tiene M(~) = 0 ⇔ M(a~) = 1 (al ser W no atómica), y como además a~ = R[W], resulta que M(R[X]) = 0 ⇔ M(R[W]) = 1. Utilizando la hipótesis inductiva se tiene que M(R[X]) = 0 ⇔ $V_m$(W) = 1. Por las tablas de verdad se concluye que M(R[X]) = 0 ⇔ $V_m$(~W) = 0 (al ser W no atómica), es decir, M(R[X]) = 0 ⇔ $V_m$(X) = 0.

Caso 6: Sea X = Y∧Z. Se tiene que R[X] = ∧, por lo que M(R[X]) = 1 ⇔ M(∧) = 1, por A∧ y nRinRd∧ se tiene M(∧) = 1 ⇔ [M(i∧) ≠ 0 y M(d∧) ≠ 0], y como además i∧ = R[Y] y d∧ = R[Z], resulta que M(R[X]) = 1 ⇔ [M(R[Y]) ≠ 0 y M(R[Z]) ≠ 0]. Utilizando la hipótesis inductiva se tiene que M(R[X]) = 1 ⇔ [$V_m$(Y) ≠ 0 y $V_m$(Z) ≠ 0]. Por la tabla de verdad se concluye que M(R[X]) = 1 ⇔ $V_m$(Y∧Z) = 1, es decir, M(R[X]) = 1 ⇔ $V_m$(X) = 1. Se concluye que M(R[X]) = 0 ⇔ $V_m$(X) = 0.

Caso 7: Sea X = Y∨Z. Se tiene que R[X] = ∨, por lo que M(R[X]) = 0 ⇔ M(∨) = 0, pero por R∨ y RiRd∨ se tiene M(∨) = 0 ⇔ [M(i∨) = 0 y M(d∨) = 0], y como además i∨ = R[Y] y d∨ = R[Z], resulta que M(R[X]) = 0 ⇔ [M(R[Y]) = 0 y M(R[Z]) = 0]. Utilizando la hipótesis inductiva se tiene que M(R[X]) = 0 ⇔ [$V_m$(Y) = 0 y $V_m$(Z) = 0]. Por la tabla de verdad se concluye que M(R[X]) = 0 ⇔ $V_m$(Y∨Z) = 0, es decir, M(R[X]) = 0 ⇔ $V_m$(X) = 0.

Caso 8: Sea X = Y→Z. Se tiene que R[X] = →, por lo que M(R[X]) = 0 ⇔ M(→) = 0, pero por R→ y nRiRd→ se tiene M(→) = 0 ⇔ [M(i→) ≠ 0 y M(d→) = 0], y como además i→ = R[Y] y d→ = R[Z], resulta que M(R[X]) = 0 ⇔ [M(R[Y]) ≠ 0 y M(R[Z]) = 0]. Utilizando la hipótesis inductiva se tiene que M(R[X]) = 0 ⇔ [$V_m$(Y) ≠ 0 y $V_m$(Z) = 0]. Por la tabla de verdad se concluye que M(R[X]) = 0 ⇔ $V_m$(Y→Z) = 0, es decir, M(R[X]) = 1 ⇔ $V_m$(X) = 1. Se concluye que M(R[X]) = 0 ⇔ $V_m$(X) = 0.

Se tiene entonces que para todos los casos M(R[X]) = 0 ⇔ $V_m$(X) = 0, quedando así probado el paso de inducción. Por el principio de Inducción se concluye que: M(R[X]) =0 ⇔ $V_m$(X) = 0.

Se ha probado entonces que para cada fórmula X de P1 y para cada función de marca *m*, existe una asignación $v_m$, tal que, M(R[X]) = 0 ⇔ $V_m$(X) = 0.

## 7.2. Proposición. *Función de marca correspondiente a una asignación* |

Para cada fórmula X de P1 y para cada asignación $v$, existe una función de marca $m_v$, tal que, $M_v(R[X]) = 0 \Leftrightarrow V(X) = 0$.

Prueba: Sea p una fórmula atómica de P1 y sea $v$ una asignación. Se define la función de marca $m_v$ del conjunto de hojas en el conjunto {0, *, 1} de la siguiente forma:

$$m_v(p)= v(p).$$

La función $m_v$ se extiende a una función $M_v$ del conjunto de nodos en el conjunto {0, *, 1}, mediante las *reglas primitivas para el forzamiento de marcas* presentadas en la sección 4. Se tiene entonces que $m_v$ es una función de marca de nodos.

Para probar que $M_v(R[X]) = 0 \Leftrightarrow V(X) = 0$, se procede por inducción sobre la *profundidad* de Ar[X].

Paso base: Supóngase que P(Ar[X]) = 0, esto significa que X es una fórmula atómica. Sea X=p, con p atómica. Como $m_v(p)= v(p)$ entonces $v(p) = 0 \Leftrightarrow m_v(p) = 0$, es decir, $M_v(R[X]) = 0 \Leftrightarrow V(X) = 0$.

Paso de inducción: Supóngase que P(Ar[X]) ≥ 1, por lo que X tiene una de las siguientes formas: ¬p, ~p, *Ip*, ¬W, ~W, Y∧Z, Y∨Z o Y→Z, donde p atómica y W no atómica.

Caso 1: Sea X = Ip, donde p es atómica. Se tiene que R[X] = I, por lo que $M_v(R[X]) = 0 \Leftrightarrow M_v(I) = 0$, pero por las reglas RI y A*a¬ se tiene $M_v(I) = 0 \Leftrightarrow M_v(p) = *$, es decir $M_v(I) = 0 \Leftrightarrow m_v(p) = *$, lo cual por definición implica $M_v(I) = 0 \Leftrightarrow v(p) = *$, utilizando la tabla de verdad resulta $M_v(I) = 0 \Leftrightarrow V(Ip) = 0$, por lo que $M_v(R[X]) = 0 \Leftrightarrow V(X) = 0$.

Caso 2: Sea X = ¬p, donde p es atómica. Se tiene que R[X] = ¬, por lo que $M_v(R[X]) = 0 \Leftrightarrow M_v(\neg) = 0$, pero por las reglas R¬ y Aa¬ se tiene $M_v(\neg) = 0 \Leftrightarrow M_v(p) = 1$, es decir $M_v(\neg) = 0 \Leftrightarrow m_v(p) = 1$, lo cual por definición implica $M_v(\neg) = 0 \Leftrightarrow v(p) = 1$, utilizando la tabla de verdad resulta $M_v(\neg) = 0 \Leftrightarrow V(\neg p) = 0$, por lo que $M_v(R[X]) = 0 \Leftrightarrow V(X) = 0$.

Caso 3: Sea X = ~p, donde p es atómica. Se tiene que R[X] = ~, por lo que $M_v(R[X]) = 0 \Leftrightarrow M_v(\sim) = 0$, pero por las reglas R~ y nRa~ se tiene $M_v(\sim) = 0 \Leftrightarrow M_v(p) \neq 0$, es decir $M_v(\sim) = 0 \Leftrightarrow m_v(p) \neq 0$, lo cual por definición implica $M_v(\sim) = 0 \Leftrightarrow v(p) \neq 0$, utilizando la tabla de verdad resulta $M_v(\sim) = 0 \Leftrightarrow V(\sim p) = 0$, por lo que $M_v(R[X]) = 0 \Leftrightarrow V(X) = 0$.

Caso 4: Sea X = ¬W, con W no atómica. Se tiene que R[X] = ¬, por lo que $M_v(R[X]) = 0 \Leftrightarrow M_v(\neg) = 0$, pero por R¬ y Aa¬ se tiene $M_v(\neg) = 0 \Leftrightarrow M_v(a\neg) = 1$, y como además a¬ = R[W], resulta que $M_v(R[X]) = 0 \Leftrightarrow M_v(R[W]) = 1$. Utilizando la hipótesis inductiva se tiene que $M_v(R[X]) = 0 \Leftrightarrow V(W) = 1$. Por la tabla de verdad se concluye que $M_v(R[X]) = 0 \Leftrightarrow V(\neg W) = 0$, es decir, $M_v(R[X]) = 0 \Leftrightarrow V(X) = 0$.

Caso 5: Sea X = ~W, con W no atómica. Se tiene que R[X] = ~, por lo que $M_v(R[X]) = 0 \Leftrightarrow M_v(\sim) = 0$, pero por R~ y nRa~ se tiene $M_v(\sim) = 0 \Leftrightarrow M_v(a\sim) \neq 0$, y como además a~ = R[W], resulta que $M_v(R[X]) = 0 \Leftrightarrow M_v(R[W]) \neq 0$. Utilizando la hipótesis inductiva se tiene que $M_v(R[X]) = 0 \Leftrightarrow V(W) \neq 0$. Por la tabla de verdad se concluye que $M_v(R[X]) = 0 \Leftrightarrow V(\sim W) = 0$, es decir, $M_v(R[X]) = 0 \Leftrightarrow V(X) = 0$.

Caso 6: Sea X de la forma Y∧Z. Se tiene que R[X] = ∧, por lo que $M_v(R[X]) = 1 \Leftrightarrow M_v(\wedge) = 1$, pero por A∧ y nRinRd∧ se tiene $M_v(\wedge) = 1 \Leftrightarrow [M_v(i\wedge) \neq 0$ y $M_v(d\wedge) \neq 0]$, y como además i∧ = R[Y] y d∧ = R[Z], resulta que $M_v(R[X]) = 1 \Leftrightarrow [M_v(R[Y]) \neq 0$ y $M_v(R[Z]) \neq 0]$. Utilizando la hipótesis inductiva se tiene que $M_v(R[X]) = 1 \Leftrightarrow [V(Y) \neq 0$ y $V(Z) \neq 0]$. Por la tabla de verdad se concluye que $M_v(R[X]) = 1 \Leftrightarrow V(Y\wedge Z) = 1$, es decir, $M_v(R[X]) = 1 \Leftrightarrow V(X) = 1$, y por lo tanto, $M_v(R[X]) = 0 \Leftrightarrow V(X) = 0$.

Caso 7: Sea X de la forma Y∨Z. Se tiene que R[X] = ∨, por lo que $M_v(R[X]) = 0 \Leftrightarrow M_v(\vee) = 0$, pero por R∨ y RiRd∨ se tiene $M_v(\vee) = 0 \Leftrightarrow [M_v(i\vee) = 0$ y $M_v(d\vee) = 0]$, y como además i∨ = R[Y] y d∨ = R[Z], resulta que $M_v(R[X]) = 0 \Leftrightarrow [M_v(R[Y]) = 0$ y $M_v(R[Z]) = 0]$. Utilizando la hipótesis inductiva se tiene que $M_v(R[X]) = 0 \Leftrightarrow [V(Y) = 0$ y $V(Z) = 0]$. Por la tabla de verdad se concluye que $M_v(R[X]) = 0 \Leftrightarrow V(Y\vee Z) = 0$, es decir, $M_v(R[X]) = 0 \Leftrightarrow V(X) = 0$.

Caso 8: Sea X de la forma Y→Z. Se tiene que R[X] = →, por lo que $M_v(R[X]) = 0 \Leftrightarrow M_v(\rightarrow) = 0$, pero por R→ y nRiRd→ se tiene $M_v(\rightarrow) = 0 \Leftrightarrow [M_v(i\rightarrow) \neq 0$ y $M_v(d\rightarrow) = 0]$, y como además i→ = R[Y] y d→ = R[Z], resulta que $M_v(R[X]) = 0 \Leftrightarrow [M_v(R[Y]) \neq 0$ y $M_v(R[Z]) = 0]$. Utilizando la hipótesis inductiva se tiene que $M_v(R[X]) = 0 \Leftrightarrow [V(Y) \neq 0$ y $V(Z) = 0]$. Por la tabla de verdad se concluye que $M_v(R[X]) = 0 \Leftrightarrow V(Y\rightarrow Z) = 0$, es decir, $M_v(R[X]) = 0 \Leftrightarrow V(X) = 0$.

Se tiene entonces que para todos los casos $M_v(R[X])=0 \Leftrightarrow V(X)=0$, quedando así probado el paso de inducción. Por el principio de Inducción se concluye que: $M_v(R[X])=0 \Leftrightarrow V(X)=0$.

Se ha probado entonces que para cada fórmula X de P1 y para cada asignación v, existe una función de marca $m_v$, tal que, $M_v(R[X])=0 \Leftrightarrow V(X)=0$.

### **7.3. Proposición.** *Caracterización semántico-deductiva de los árboles de forzamiento semántico trivalentes* |

1. La fórmula X es válida desde el punto de vista de los árboles si y solamente si X es válida desde el punto de vista de la semántica tri-valuada.

2. La fórmula X es válida desde el punto de vista de los árboles si y solamente si X es un teorema de los sistemas deductivos paraconsistentes LPcAt y P1.

Prueba: Supóngase que X no es válida desde el punto de vista de los árboles, entonces existe una función de marca $m$, tal que $M(R[X]) = 0$. Se tiene entonces, por la proposición 7.1, que existe una asignación $v_m$, tal que $V_m(X) = 0$, y por lo tanto, X no puede ser válida desde el punto de vista de la semántica trivaluada.

Supóngase ahora que X no es válida desde el punto de vista de la semántica tri-valuada, entonces existe una asignación v, tal que $V(X) = 0$. Se tiene entonces, por la proposición 7.2, que existe una función de marca $m_v$, tal que $M_v(R[X]) = 0$, y por lo tanto, X no puede ser válida desde el punto de vista de los árboles.

Se concluye así que, X es válida desde el punto de vista de los árboles si y solamente si X es válida desde el punto de vista de la semántica tri-valuada.

De (Sette, 1973; Sierra, 2003), se sabe que los teoremas de LPcAT y P1 son exactamente las fórmulas validas desde el punto de vista de la semántica tri-valuada, entonces, como consecuencia de la parte 1 se tiene que, la fórmula X es válida desde el punto de vista de los árboles si y solamente si X es un teorema de los sistemas deductivos LPcAt y P1.

## 8. ILUSTRACIONES |

Para las marcas de los nodos se siguen las convenciones presentadas en la sección 5.14 (presentación gráfica de las reglas). Para el análisis de validez se utilizan los métodos señalados en la sección 4.7 (validez).

### 8.1. Ilustración. *Trivialización* |

En la figura 21 se muestra un árbol de forzamiento *mal marcado* para la fórmula A-válida IA→[¬A→(A→B)], donde A es una fórmula atómica, así como el árbol de forzamiento *bien marcado* de la fórmula A-inválida ¬A→(A→B).

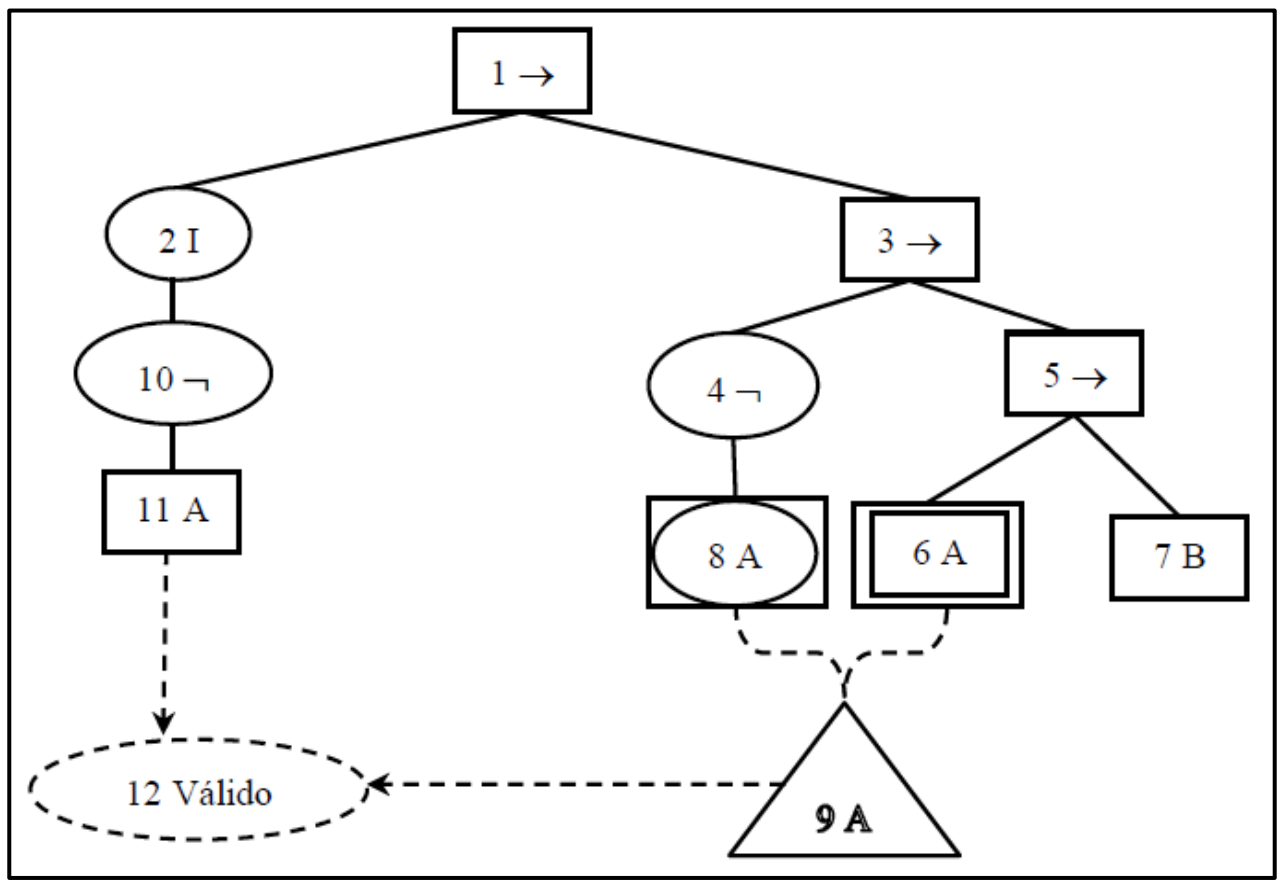

Figura 21

Justificaciones:

1. RR — 2, 3. R→ en 1 — 4, 5. R→ en 3

6, 7. R→ en 5 — 8. A¬ en 4

9. nRnAat en 8 y 6 — 10. IM en 4

11. A¬AI en 10 y 2 — 12. RR-DM en 1, 9 y 11

Observar, además, que los pasos de 3 a 9 determinan un árbol *bien marcado* para la fórmula ¬A→(A→B), y en consecuencia ésta es A-inválida, además, las marcas de las hojas indican que la fórmula es refutada por la siguiente asignación: v(A)=*, v(B)=0. En resumen, cuando se tiene una fórmula atómica y su negación débil, el sistema P1 colapsa (se trivializa, todas las fórmulas se derivan de una contradicción) si la fórmula atómica no puede ser aceptada por defecto (es incompatible con su negación).

### 8.2. Ilustración. *Reducción al absurdo* |

En la figura 22 se muestra un árbol de forzamiento *mal marcado* para la fórmula A-válida IB→{(A→B)→[(A→¬B)→¬A]}, donde B es una fórmula atómica, así como el árbol de forzamiento *bien marcado* de la fórmula A-inválida (A→B)→[(A→¬B)→¬A].

Justificaciones:

1. RR — 2, 3. R→ en 1

4, 5. R→ en 3 — 6, 7. R→ en 5

8. R¬ en 7 — 9. IM en 8

10. nRiA→ en 9 y 4 — 11. IM en 8

12. nRiA→ en 11 y 6 — 13. A¬ en 12

14. nRnAat en 10 y 13 — 15. IM en 12

16. A¬AI en 15 y 2 — 17. RR-DM en 1, 16 y 14

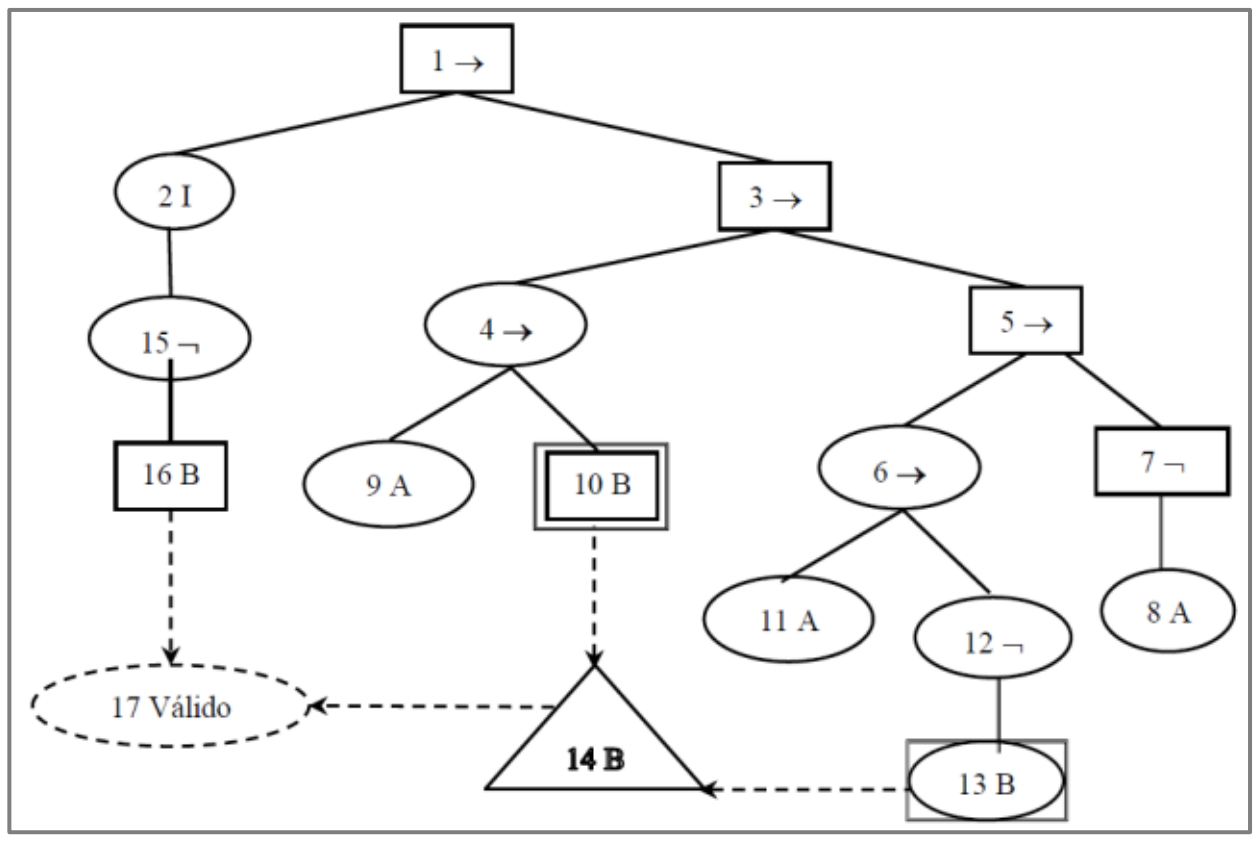

Figura 22

Observar que los pasos de 3 a 14 determinan un árbol *bien marcado* para la fórmula (A→B)→[(A→¬B)→¬A], y en consecuencia ésta es A-inválida, y además, las marcas de las hojas indican que la fórmula es refutada por la siguiente asignación: v(B) = * y v(A) = 1.

El análisis realizado a la fórmula IB→{(A→B)→[(A→¬B)→¬A]}, corresponde a la formalización del siguiente argumento: Si una propuesta genera consecuencias que ni son aceptadas (son cuestionadas) (A→¬B) ni son rechazadas (A→B), entonces no se puede asegurar que la propuesta no sea aceptada (sea cuestionada) (¬A), salvo que sea forzosa la aceptación o el rechazo de las consecuencias (IB).

## **8.3. Ilustración.** *Negación de la conjunción* |

En la figura 23 se muestra un árbol de forzamiento con la raíz *forzosamente marcada con 1* para la fórmula A-válida ¬(A∧B)→(¬A∨¬B), donde A y B son fórmulas atómicas.

Justificaciones:

1. ORd
2, 3. R∨ en 1
4. R¬ en 2
5. R¬ en 3
6. IM en 4
7. IM en 5
8. nRinRd∧ en 6 y 7
9. Aa¬ en 8
10. ORd-Ri→ en 1, 9
11. RM1 en 10

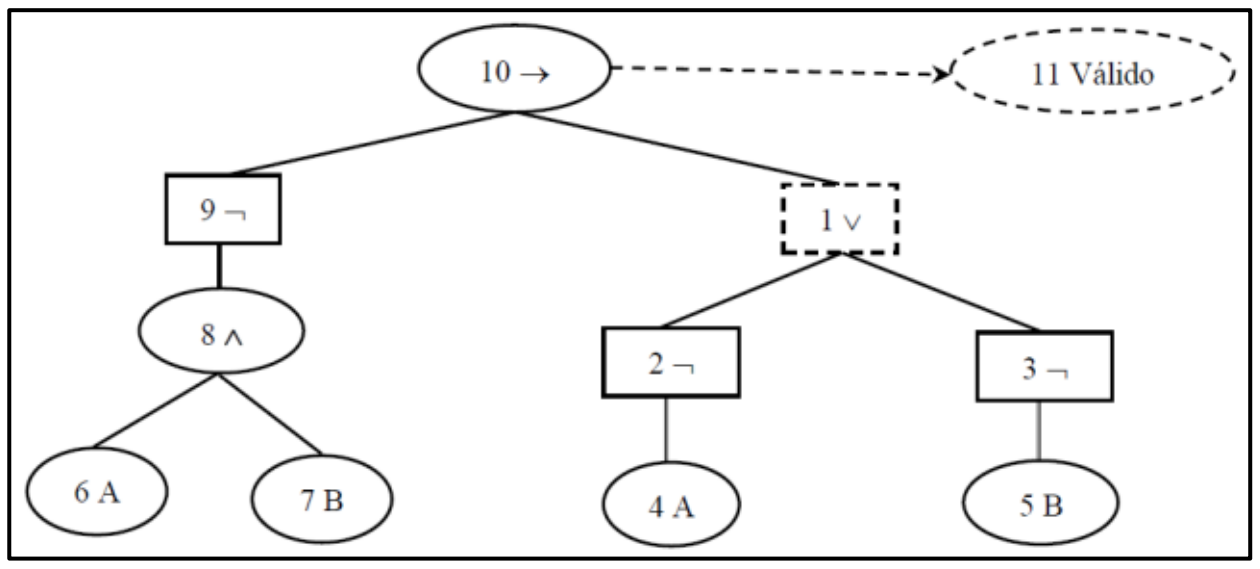

Figura 23

El análisis realizado a la fórmula ¬(A∧B)→(¬A∨¬B), corresponde a la formalización del siguiente argumento: Si no se acepta la ejecución de ambos proyectos (¬(A∧B)), entonces la ejecución de al menos uno de los proyectos no es aceptada (¬A∨¬B).

## **8.4. Ilustración.** *Disyunción de negaciones* |

En la figura 24 se muestra un árbol de forzamiento *mal marcado* para la fórmula A-válida (IA∧IB)→[(¬A∨¬B)→¬(A∧B)], donde A y B son fórmulas atómicas, así como el árbol de forzamiento *bien marcado* de la fórmula A-inválida (¬A∨¬B)→¬(A∧B).

Justificaciones:

1. RR
2, 3. R→en 1
4, 5. R→ en 3
6. R¬ en 5
7, 8. A∧ en 6
9, 10. A∧ en 2
11. nRa¬AI en 7 y 9
12. IM en 11
13. Aa¬ en 12
14. RiA∨ en 13 y 4
15. A¬ en 14
16. nRnAat en 8 y 15
17. IM en 14
18. A¬AI en 17 y 10
19. RR-DM en 1, 16 y 18

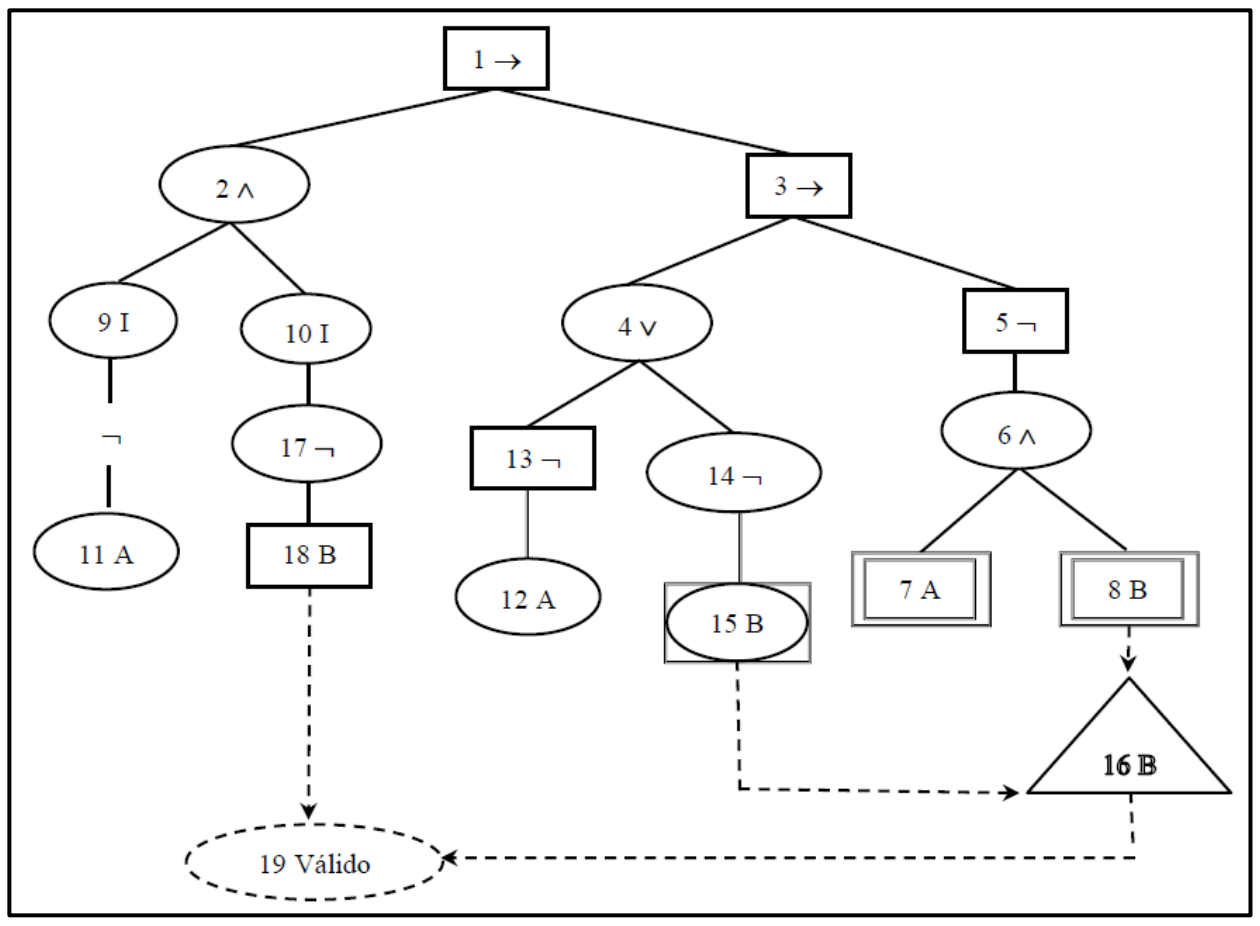

Figura 24

Observar, que los pasos de 3 a 8, junto con 12 a 16 (justificando 12 como opción de marca OM), determinan un árbol *bien marcado* para la fórmula (¬A∨¬B)→¬(A∧B), y en consecuencia ésta es A-inválida, y además, las marcas de las hojas indican que la fórmula es refutada por la siguiente asignación: v(A)=1 y v(B)=*.

El análisis realizado a la fórmula (IA∧IB)→[(¬A∨¬B)→¬(A∧B)], corresponde a la formalización del siguiente argumento: Si de al menos uno de los productos no se acepta la compra (¬A∨¬B), entonces no se puede asegurar que no se acepte comprar ambos productos (¬(A∧B)). Esto se debe a que acepta la compra de A y aunque no acepta la compra de B, tampoco la rechaza. Pero, se puede asegurar que no se acepte comprar ambos productos si es forzosa la aceptación o el rechazo de la compra de cada uno de los productos (IA∧IB).

## **8.5. Ilustración.** *No transmisión de la incompatibilidad* |

En la figura 25 se muestra un árbol de forzamiento *bien marcado* para la fórmula A-inválida ((A∧IA)∧(A→B))→IB, donde A y B son fórmulas atómicas.

Justificaciones:

1. RR
2, 3. R→en 1
4. RI en 3
5. A*a¬ en 4
6, 7. A∧ en 2
8, 9. A∧ en 6
10. nRa¬AI en 8, 9
11. IM en 10
12. nRiA→ en 11, 7
13. Aa¬ en 10
14. ABM

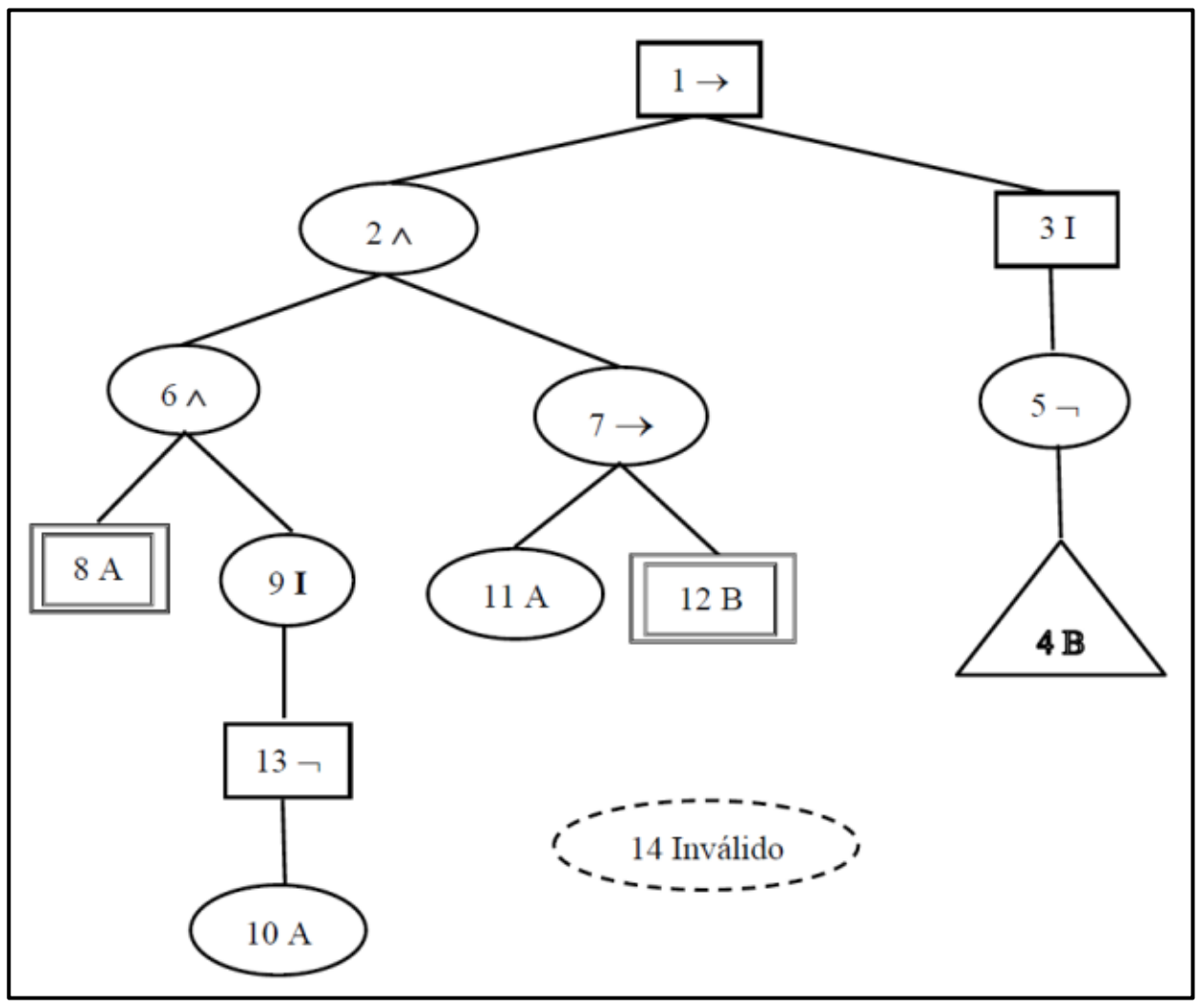

Figura 25

Observar que las marcas de las hojas indican que la fórmula es refutada por la siguiente asignación: v(A)=1 y v(B)=*.

En resumen, cuando una fórmula atómica es incompatible con su negación, no se puede asegurar que sus consecuencias lo sean.

## 8.6. Ilustración. *Tercero excluido* |

En la figura 26 se muestra un árbol de forzamiento *con la raíz marcada con 1* para la fórmula A-válida A∨¬A (tercero excluido), donde A es una fórmula arbitraria.

Justificaciones:

1. ORi — 2. IM en 1 — 3. Ra¬ en 2

4. ORi-Ad∨ en 1 y 3 — 5. RM1 en 4

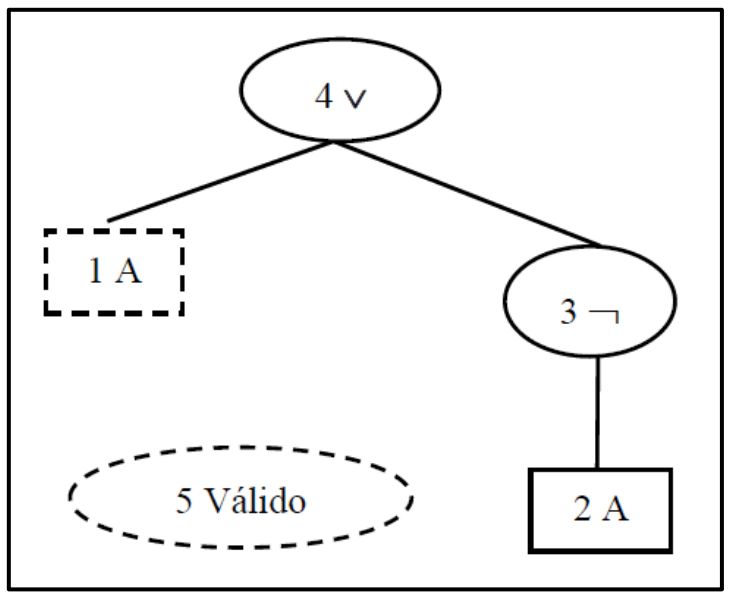

Figura 26

Resulta entonces que, en P1 vale el principio del tercero excluido.

## 8.7. Ilustración. *Trivialización con negación débil* |

En la figura 27 se muestra un árbol de forzamiento *con la raíz marcada con 1* para la fórmula A-válida IA→(A→(¬A→B)), donde A es una fórmula atómica y B es una fórmula arbitraria, así como el árbol de forzamiento *bien marcado* de la fórmula A-inválida A→(¬A→B). En este caso, la prueba de validez se realiza utilizando *forzamiento directo*, mientras que en la ilustración 8.1, la prueba de validez se realiza utilizando *forzamiento indirecto*.

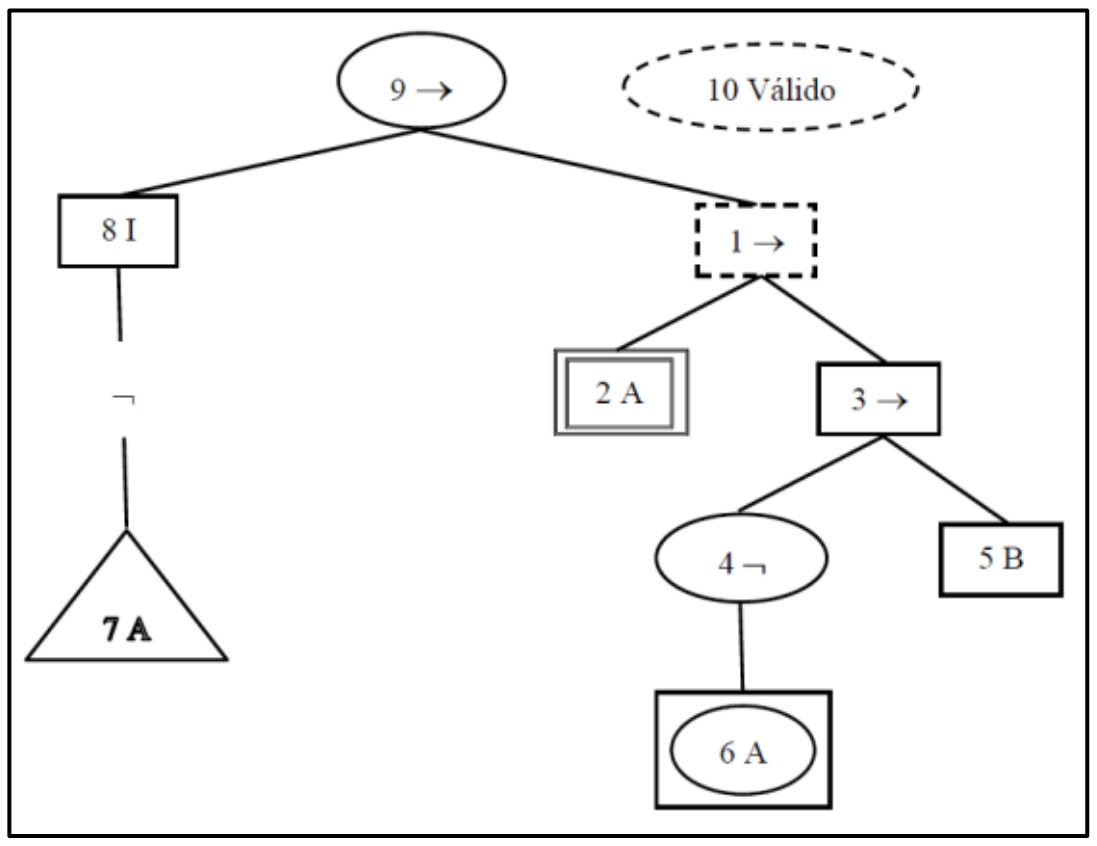

Figura 27

Justificaciones:

1. ORd — 2, 3. R→ en 1 — 4, 5. R→ en 3

6. A¬ en 4 — 7. nAnRat en 2, 6 — 8. A*a¬I en 7

9. ORd-Ri→ en 1, 8 — 10. RM1 en 9

Observar que los pasos de 1 a 7, determinan un árbol *bien marcado* para la fórmula A→(¬A→B), y en consecuencia ésta es A-inválida, y además, las marcas de las hojas indican que la fórmula es refutada por la siguiente asignación: v(A)=* y v(B)=0.

Observar que, si A no es atómica entonces existiría una contradicción entre los pasos 2 y 6, resultando que esta última fórmula es válida. En resumen, en la figura 27 se verifica que el sistema P1 es paraconsistente a nivel atómico respecto al cuestionamiento (negación débil ¬).

## 8.8. Ilustración. *Trivialización con negación fuerte* |

En la figura 28 se muestra un árbol de forzamiento *con la raíz marcada con 1* para la fórmula A-válida A→(~A→B), donde A y B son fórmulas arbitrarias. Resultando que P1 no es paraconsistente respecto a la negación fuerte ~.

Justificaciones:

1. RR — 2, 3. R→ en 1 — 4, 5. R→ en 3

6. A~ en 4 7. RR-DM en 1, 2 y 6.

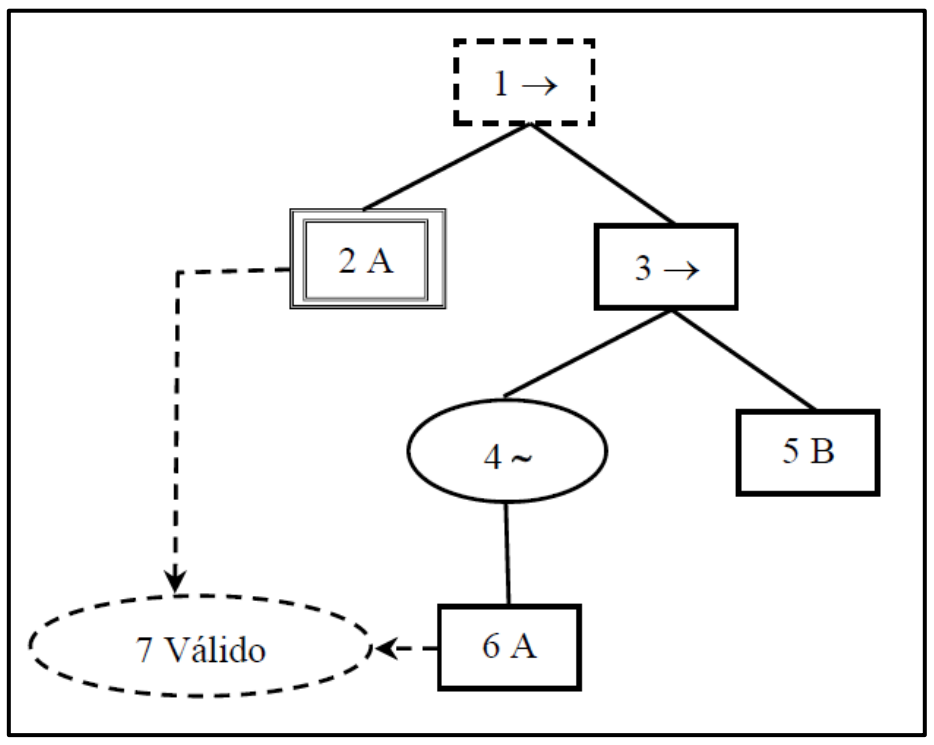

Figura 28

## 9. CONCLUSIONES |

Con los árboles de forzamiento semántico trivalentes para el sistema P1, la validez de una fórmula se puede determinar de manera visual y completamente mecánica, por ejemplo, mediante un algoritmo se recorre el árbol de la fórmula, buscando en cada nodo la aplicación de una regla para marcarlo. Cuando la fórmula es inválida, es decir, cuando el árbol de la fórmula está bien marcado, entonces la lectura de las marcas de las hojas (fórmulas atómicas), proporciona la asignación de valores de verdad, con la cual se refuta la validez de la fórmula.

Desde el punto de vista didáctico, la simplicidad de las reglas para el forzamiento de marcas, en comparación con las deducciones de los sistemas axiomáticos P1 y LPcAt, hacen de los árboles de forzamiento semántico trivalentes, una herramienta de trabajo muy útil.

## 9. REFERENCIAS